\documentclass[a4paper,12pt,reqno]{amsart}
\usepackage{amssymb}
\usepackage[dvips]{graphicx}

\usepackage{hyperref}
\usepackage[all]{xy}

\input epsf.tex
\newdimen\xsize
\newdimen\oldbaselineskip
\newdimen\oldlineskiplimit
\def\restorelineskip{\baselineskip=\oldbaselineskip%
\lineskiplimit=\oldlineskiplimit}
\def\putm[#1][#2]#3{
\hbox{\vbox to 0pt{\parindent=0pt%
\vskip#2\xsize\hbox to0pt{\hskip#1\xsize $#3$\hss}\vss}}}%
\long\def\Line#1{\hbox to \hsize{#1}}
\def\putt[#1][#2]#3{
\vbox to 0pt{\noindent\hskip#1\xsize\lower#2\xsize%
\vtop{\restorelineskip#3}\vss}}

\makeatletter
\def\xbig[#1]#2{{\hbox{$\m@th\left#2\vbox to#1\xsize{}%
\right.\n@space$}}}
\makeatother
\def\xlar[#1]#2{%
\smash{\mathop{ \hbox to #1\xsize{\leftarrowfill}}\limits^{#2}}}
\def\xrar[#1]#2{%
\smash{\mathop{ \hbox to #1\xsize{\rightarrowfill}}\limits^{#2}}}
\def\xline[#1]{\hbox to #1\xsize{\leaders\hrule\hfill}}

\DeclareFontFamily{U}{rsf}{\skewchar\font'177}%
\DeclareFontShape{U}{rsf}{m}{n}{<-6>rsfs5<6-8>rsfs7<8->rsfs10}{}%
\DeclareFontShape{U}{rsf}{b}{n}{<-6>rsfs5<6-8>rsfs7<8->rsfs10}{}%
\DeclareMathAlphabet\RSFS{U}{rsf}{m}{n}
\SetMathAlphabet\RSFS{bold}{U}{rsf}{b}{n}
\def\sf#1{{\mathsf{#1}}}

\def\slsf{\slshape \sffamily }

\def\msmall#1{\mathchoice{\hbox{\small$\displaystyle {#1}$}}{#1}{#1}{#1}}

\def\bb{{\mathbb B}}

\def\cc{{\mathbb C}}

\def\sss{{\mathbb S}}
\def\nn{{\mathbb N}}

\def\pp{{\mathbb P}}

\def\loc{{\sf{loc}}}

\def\area{\sf{area}}

\def\const{\sf{const}}

\def\codim{\sf{codim}\,}

\def\deg{\sf{deg}\,}

\def\dim{\sf{dim}\,}

\def\ma{\sf{MA}}
\def\sym{\sf{Sym}}

\def\lim{\mathop{\sf{lim}}}

\def\max{\sf{max}}
\def\min{\sf{min}}

\def\pr{\sf{pr}}

\def\Sing{\sf{Sing}\,}

\def\sup{\sf{sup}\,}

\def\vol{\sf{Vol}}

\def\eps{\varepsilon}

\def\<{\langle}\let\la=\<
\def\>{\rangle}\let\ra=\>
  
\def\comp{\Subset}
\def\d{\partial}

\def\ddef{\mathrel{{=}\raise0.3pt\hbox{:}}}
\def\deff{\mathrel{\raise0.3pt\hbox{\rm:}{=}}}

\def\fraction#1/#2{\mathchoice{{\msmall{ #1\over#2}}}%
{{ #1\over #2 }}{{#1/#2}}{{#1/#2}}}
\def\norm#1{\left\Vert{#1}\right\Vert}

\def\le{\leqslant}

\def\emptyset{\varnothing}

\def\longpoints{\leaders\hbox to 0.5em{\hss.\hss}\hfill \hskip0pt}
\def\stateskip{\smallskip}
\def\state#1. {\stateskip\noindent{\bf#1. }} 
\def\statep#1. {\stateskip\noindent{\bf#1 }} 
\def\proof{\state Proof. }

\def\Chi{\raise 2pt\hbox{$\chi$}}

\def\ie{\hskip1pt plus1pt{\sl i.e.\/,\ \hskip1pt plus1pt}}

\def\sli{{\sl i)} } 
\def\slii{{\sl i$\!$i)} } 
\def\sliii{{\sl i$\!$i$\!$i)} }

\def\Chi{\raise 2pt\hbox{$\chi$}}

\let\phI=\phi\let\phi=\varphi\let\varphi=\phI
\let\cal=\mathcal

\def\calc{{\cal C}}

\def\calf{{\cal F}}

\def\calm{{\cal M}}
\def\caln{{\cal N}}
\def\calo{{\cal O}}

\def\calu{{\cal U}}

\def\eps{\varepsilon}

\def\comp{\Subset}
\def\d{\partial}

\def\1{{1\mkern-5mu{\rom l}}}

\def\ge{\geqslant}

\def\fraction#1/#2{\mathchoice{{\msmall{ #1\over#2}}}%
{{ #1\over #2 }}{{#1/#2}}{{#1/#2}}}

\def\le{\leqslant}

\def\emptyset{\varnothing}

\def\qed{\ \ \hfill\hbox to .1pt{}\hfill\hbox to .1pt{}\hfill $\square$\par}

\def\comment#1\endcomment{}

 \belowdisplayskip=\abovedisplayskip

\def\lineeqqno(#1){\hfill\llap{\vbox to 10pt%
{\vss\begin{align} \eqqno(#1)\end{align}\vss}}\vskip1pt}

\advance\headheight 1.2pt

\def\ShowwLLabel#1{}

\def\thechpt{\Roman{chpt}}

\def\newchapt[#1]#2{\newpage%
\refstepcounter{chpt}\setcounter{subsection}{0}%
\setcounter{thm}{0}\setcounter{defi}{0}%
\setcounter{rema}{0}\setcounter{exrc}{0}%
\renewcommand{\thesubsection}{\thechpt.\arabic{subsection}}%
\section*{\begin{center}\huge \bf Chapter \thechpt\\
#2 \end{center}}\label{#1}%
\ \smallskip%
\markboth{Chapter \thechpt}{#2}%
}
\def\newsect[#1]#2{\refstepcounter{section}\setcounter{equation}{0}%
\renewcommand{\thesubsection}{\arabic{section}.\arabic{subsection}}%
\section*{\arabic{section}.
#2}\vspace{-20pt}\label{#1}\vspace{20pt}%
\markboth{Section \arabic{section}}{#2}}

\def\newlect[#1]#2{\refstepcounter{section}%
\renewcommand{\thesubsection}{\arabic{section}.\arabic{subsection}}%
\section*{Lecture \arabic{section}\\
#2}\label{#1}%
\markboth{Lecture \arabic{section}}{#2}}

\def\newprg[#1]#2{\refstepcounter{subsection}%
\subsection*{{\thesubsection.\ #2}} \label{#1}%
}

\def\newappx[#1]#2{%
\refstepcounter{appx}\setcounter{section}{0}%
\renewcommand{\thesubsection}{A\arabic{appx}.\arabic{subsection}}%
\section*{Appendix \arabic{appx}\\ #2}
\label{#1}%
\markboth{Appendix A\arabic{appx}}{#2}
}

\newtheorem{thm}{Theorem}[section]
   \def\newthm#1{\begin{thm}\label{#1}}

\newtheorem{nnthm}{Theorem}
   \def\newthm#1{\begin{nnthm}\label{#1}}

\newtheorem{lem}{Lemma}[section]
   \def\newlemma#1{\begin{lem} \label{#1}}

\newtheorem{prop}{Proposition}[section]
   \def\newprop#1{\begin{prop}\label{#1}}

\newtheorem{nnprop}{Proposition}
   \def\newprop#1{\begin{nnprop}\label{#1}}

\newtheorem{corol}{Corollary}[section]
   \def\newcorol#1{\begin{corol} \label{#1}}

\newtheorem{nncorol}{Corollary}
   \def\newcorol#1{\begin{nncorol} \label{#1}}

\newtheorem{nndefi}{Definition}
   \def\newdefi#1{\begin{nndefi} \label{#1}\rm }

\newtheorem{defi}{Definition}[section]
   \def\newdefi#1{\begin{defi} \label{#1}\rm }

\newtheorem{exmp}{Example}[section]
   \def\newexmp#1{\begin{exmp} \label{#1}\rm }

\newtheorem{exrc}{Exercise}
   \def\newexrc#1{\begin{exrc} \label{#1}\rm }

\newtheorem{rema}{Remark}[section]
   \def\newrema#1{\begin{rema} \label{#1}\rm }

\newtheorem{nnrema}{Remark}
   \def\newrema#1{\begin{nnrema}\label{#1}}

\def\eqqno(#1){\label{(#1)}}
\def\eqqref(#1){(\ref{(#1)})}

\title{Weak normality of families of meromorphic\\
mappings and bubbling in higher dimensions}
\author {S. Ivashkovich, F. Neji}
\date{\today}

\address{
Universit\'e de Lille-1, UFR de Math\'ematiques, 59655 Villeneuve
d'Ascq, France} \email{ivachkov@math.univ-lille1.fr}
\address{
Universit\'e de Lille-1, UFR de Math\'ematiques, 59655 Villeneuve
d'Ascq, France} \email{Fethi.neji@math.univ-lille1.fr}
\begin{document}

\begin{abstract}
Our primary goal in this paper is to understand wether the  sets of
normality of families of meromorphic mappings between general
complex manifolds are pseudoconvex or not. It turns out that the
answer crucially depends on the type of convergence one is
interested in. We examine three natural types of convergence
introduced by one of us earlier and prove pseudoconvexity of sets of
normality for a large class of target manifolds for the so called
{\slsf weak} and {\slsf gamma} convergencies. Furthermore we determine
the structure of the exceptional components of the limit of a weakly/gamma
but not strongly converging sequence, they turn to be {\slsf rationally
connected}. This observation allows to determine effectively when a
weakly/gamma converging sequence fails to converge strongly. An application
to the Fatou sets of meromorphic self-maps of compact complex surfaces
is given.
\end{abstract}

\maketitle

\setcounter{tocdepth}{1}
\tableofcontents

\newsect[INTR]{Introduction}

\newprg[INTR.conv]{Convergence of meromorphic mappings}
When one works with sequences of meromorphic functions and, more generally,
mappings one finds himself bounded to consider several notions of
their convergence. Some of these notions were introduced in
\cite{Fu} and \cite{Iv2}, we shall recall the essentials below.
An important question is: what can be said
about the maximal open sets where the given sequence converge?
It occurs that pseudoconvexity or not of domains of
convergence/normality in the case of meromorphic mappings crucially
depends on the type of convergence one is looking for.

\smallskip Now let briefly describe the ways one can define what does
it means that a sequence $\{f_k\}$ of meromorphic mappings between
complex manifolds $U$ and $X$ converges. We start with the most obvious one.
A sequence $\{f_k\}$ of
meromorphic mappings between complex manifolds $U$ and $X$ is said
to converge {\slsf strongly} to a meromorphic map $f$ if the graphs
$\Gamma_{f_k}$ converge  over compacts in $U$ to the graph
$\Gamma_f$ in Hausdorff metric. Our first result shows that
$\Gamma_{f_k}$ then converge
to $\Gamma_f$ in a stronger {\slsf topology of cycles}.

\begin{nnthm}
\label{str-cycl}
If $f_k$ strongly converge to $f$ then for every compact $K\comp U$
the volumes $\Gamma_{f_k}\cap (K\times X)$ are uniformly bounded and therefore
$\Gamma_{f_k}$ converge to $\Gamma_f$  in the topology of cycles.
\end{nnthm}

This type of convergence is natural and has some nice features. For
example the strong limit of a sequence of holomorphic maps is
holomorphic and vice versa, if the limit $f$ is holomorphic then for
every compact $K\comp U$ all $f_k$ for $k\gg 1$ are holomorphic in a
neighborhood of $K$ and uniformly converge there to $f$. This
statement was called the Rouch\'e Principle in \cite{Iv2}.

\smallskip But strong convergence has also some disadvantages. The first, crucial
for us is the fact that the sets of strong convergence, \ie maximal
open subsets of $U$ where a given sequence converges strongly on
compacts, are not pseudoconvex in general. Moreover, the sets of
strong normality (see later on) of families of meromorphic mappings
can be just arbitrary, see Example \ref{norm-exmp}. Also if one
takes $X=\pp^N$ the "most immediate" notion of convergence doesn't
correspond to the strong one.

\smallskip Therefore in \cite{Iv2} along with the notion of strong
convergence we proposed two weaker ones. We say that $f_k$ converge
{\slsf weakly} to $f$ if they converge strongly to $f$ on compacts
outside of some analytic set $A$ in $U$ of codimension at least two.
It turns out that this $A$ can be taken to be the set $I_f$ of
points of indeterminacy of the limit map $f$ and then for every
compact $K$ in $U\setminus I_f$ all weakly converging to $f$
mappings $f_k$ will be holomorphic on $K$ (for $k$ big enough) and
converge to $f$ uniformly on $K$, see Remark \ref{weak-rem}.

\smallskip One more notion of convergence from \cite{Iv2}, which we
need to recall here, is the gamma convergence
($\Gamma$-convergence). We say that $f_k$ {\slsf gamma}-converge to
$f$ if they strongly converge to $f$ outside of an analytic set (now
it can be of codimension one) and for every divisor $H$ in $X$ and
every compact $K\comp U$ the intersections $f_k^*H\cap K$ have
bounded volume counted with multiplicities, see more about the last
condition in section \ref{PSCN.gamma}.

\begin{nnrema} \rm
Strong convergence (or {\slsf s} - convergence) will be denoted by
$f_k\to f$, the weak one (or {\slsf w} - convergence)  as
$f_k\rightharpoonup f$, and $\Gamma$ - convergence as
$f_k\xrightarrow{\Gamma}f$. Note that in the second and third
definitions we suppose that the limit $f$ is defined and meromorphic
on the whole of $U$ if, even, the convergence takes place only  on
some part of $U$. In the first case the limit exists on the whole of
$U$ automatically.
\end{nnrema}

For the better understanding of these notions let us give a
description of the listed types of convergence in the case when $X$
is projective, \ie imbeds into $\pp^N$ for some $N$. In that special
case the notions of convergence listed above permit an explicit
analytic description as follows. Every meromorphic mapping $f$ with
values in $\pp^N$ can be locally represented by an $(N+1)$-tuple of
holomorphic functions

\begin{equation}
\eqqno(pn-pres) f(z) = [f^0(z):...:f^N(z)],
\end{equation}
where not all of $f^0,...,f^N$ are identically zero, see section \ref{PN}. More
precisely, if $f:U\to \pp^N$ is a meromorphic mapping then for every
point $x_0\in U$ there exists a neighborhood $V\ni x_0$ and
holomorphic functions $f^0,...,f^N$ in $V$ satisfying
\eqqref(pn-pres). If the zero sets of $f^j$ contain a common divisor
then we can divide all $f^j$ by its equation and get a
representation such  that $GCD(f^1,...,f^N)=1$ in every $\calo_x$,
$x\in V$. In that case the indeterminacy set of $f$ is

\begin{equation}
\eqqno(codim) I_f\cap V = \{z\in V: f^0(z)=...=f^N(z) = 0\}
\end{equation}
and has codimension at least two. Representation \eqqref(pn-pres)
satisfying \eqqref(codim) is  called {\slsf reduced}. We shall prove the following

\begin{nnthm}
\label{conv-pn} Let $\{f_k\}$ be a sequence of meromorphic mappings
from a complex manifold $U$ to $\pp^N$. Then:

\smallskip\sli $f_k\xrightarrow{\Gamma} f$ if and only if for any point
$x_0\in U$ there exists a neighborhood $V\ni x_0$, reduced

representations $f_k = [f_k^0:...:f_k^N]$ and {\slsf not
necessarily reduced} representation $f = [f^0:$

$...:f^N]$ such that for every $0\le j\le N$ the sequence
$f_k^j$ converges to $f^j$ uniformly on $V$;

\smallskip\slii $f_k\rightharpoonup f$ if and only if
$f_k\xrightarrow{\Gamma} f$ and the limit representation $f =
[f^0:...:f^N]$ is {\slsf reduced};

\smallskip\sliii $f_k\to f$ if and only if $f_k\rightharpoonup f$ and
corresponding {\slsf non-pluripolar}  Monge-Amp\`ere masses

converge, \ie for every $1\le p\le n=\dim U$ one has

\begin{equation}
\eqqno(monge1) \left(dd^c\norm{z}^2\right)^{n-p}\wedge
\left(dd^c\ln\norm{f_k}^2\right)^p \to
\left(dd^c\norm{z}^2\right)^{n-p}\wedge
\left(dd^c\ln\norm{f}^2\right)^p
\end{equation}

weakly on compacts in $U$.
\end{nnthm}
Here in  \eqqref(monge1) we suppose that $V=\Delta^n$, $z_1,...,z_n$
are standard coordinates and $\norm{f}^2= |f^0|^2 + ... +|f^N|^2$, \ie
$dd^c\ln\norm{f}^2$ is the pullback of the Fubini-Study form by $f$.
Non-pluripolar MA mass of $\ln\norm{f}^2$ of order $p$ in $V$ here means
\begin{equation}
\eqqno(ma-non-pol)
\int\limits_{V\setminus I_f}\left(dd^c\norm{z}^2\right)^{n-p}\wedge
\left( dd^c\ln\norm{f}^2\right)^p,
\end{equation}
where $I_f$ is given by \eqqref(codim), \ie is the indeterminacy set of $f$.

\begin{nnrema} {\bf a)} \rm Reducibility or not of the limit representation
$f = [f^0:...:f^N]$ in this theorem doesn't depend on the choice of
{\slsf converging} representations $f_k = [f_k^0:...:f_k^N]$,
provided they are taken to be reduced (the last can be
assumed always). Indeed, any other reduced
representation of $f_k$ has the form $f_k =
[g_kf_k^0:...:g_kf_k^N]$, where $g_k$ are holomorphic and nowhere
zero. If the newly chosen representations converge to some
representation of $f$ then $g_k$ must converge, say to $g$, and this
$g$ is nowhere zero by Rouch\'e's theorem. Therefore the obtained
representation of the limit is $f = [gf^0:...:gf^N]$ and it is
reduced if and only if $f = [f^0:...:f^N]$ was reduced.

\smallskip\noindent{\bf b)} The case when the representation
$f = [f^0:...:f^N]$ of the limit is not necessarily reduced was
studied for mappings with values in $\pp^N$ by H. Fujimoto in
\cite{Fu}, who called it {\slsf meromorphic}, or {\slsf
m}-convergence. According to the part (\sli of our theorem it turns
out that our $\Gamma$-convergence (in the case of $X=\pp^N$) is
equivalent to {\slsf m}-convergence of Fujimoto.
\end{nnrema}

\newprg[INT.norm]{Sets of normality and Bloch-Montel type criterion}

\smallskip In this paper we consider two classes of complex manifolds:
projective and Gauduchon, the last is the class of complex manifolds
carrying a $dd^c$-closed metric form - a Gauduchon form.
Let $\calf $ be a family of meromorphic mappings between complex
manifolds $U$ and $X$. $\calf $ is said to be {\slsf
strongly/weakly} or {\slsf gamma}  normal if from every sequence of
elements of $\calf$ one can extract a subsequence converging on
compacts in $U$ in the corresponding sense. The maximal open subset
$\caln_{\calf}\subset U$ on which $\calf$ is normal is called the
set of normality. As it was already told the sets of strong
normality could be arbitrary.  In subsection \ref{PSCN.w-conv} we
prove the following

\begin{nnthm}
\label{propa} Let $U$ be a domain in a Stein manifold $\hat U$ such
that $\hat U$ is an envelope of holomorphy of $U$ and let $f_k:\hat
U\to X$ be a weakly converging on $U$ sequence of meromorphic
mappings with values in a disk-convex complex manifold $X$. Then:

\smallskip (a) If the weak limit $f$ on $f_k$ meromorphically
extends from $U$ to $\hat U$ then  $f_k$  weakly

\quad converge to $f$ on the whole of $\hat U$.

\smallskip (b) If, in addition, the manifold $X$ carries
a pluriclosed metric form then  the weak  limit

\quad $f$ of $f_k$ meromorphically extends to $\hat U$ and then the
part (a) applies.
\end{nnthm}

As a result the sets of weak normality are locally pseudoconvex
provided the target is disk-convex and Gauduchon. Recall that an open
subset $\caln$ of a complex manifold $U$ is called locally pseudoconvex
if for every point $p\in \d \caln$ there exists a Stein neighborhood $V$ of
$p$ in $U$ such that $V\cap \caln$ is Stein.

\begin{nncorol}
\label{weak-norm} Let $\calf \subset \calm (U,X)$ be a family of
meromorphic mappings from a complex manifold $U$ to a disk-convex
Gauduchon manifold $X$. Then the set of weak normality
$\caln_{\calf}$ of $\calf$ is locally pseudoconvex. If $\calf =\{f_k\}$ is a
sequence then the set of its weak convergence is locally pseudoconvex.
\end{nncorol}

\begin{nnrema} \rm
This corollary clearly follows from Theorem
\ref{propa}. Sets of $\Gamma$-normality are also locally
pseudoconvex under the same assumptions, see Proposition
\ref{gamma-norm} in section \ref{MER-MAP}.
\end{nnrema}

As one more supporting argument in favor of weak convergence we
prove in section \ref{NORM} the following normality criterion.

\begin{nnthm}
\label{montel} Let $\{H_i\}_{i=0}^{d}$, $d\ge 1$, be hypersurfaces in projective
manifold $X$ such that $Y\deff X\setminus \bigcup_{i=0}^dH_i$ is hyperbolically
imbedded to $X$. Let $\calf $ be a family of
meromorphic mappings from a complex manifold $U$ to $X$ such that:

\smallskip\sli  for every $i=0,...,d$ and every compact $K\comp U$ the
volumes $f^*H_i\cap K$ counted with

multiplicities are uniformly bounded for $f\in \calf$;

\smallskip\slii $\calf$ uniformly separates every pair $H_i, H_j$, $0\le i
< j \le d$.

\smallskip\noindent Then the family $\calf$ is weakly normal on $U$.
\end{nnthm}

Conditions (\sli and (\slii are explained in section \ref{NORM},
they are intuitively clear and more or less necessary. The
classically known case of a system of divisors with hyperbolically
imbedded complement is $2N+1$ hypersurfaces in $\pp^N$ in general
position - Theorem of Bloch, see \cite{Gr}. A criterion for {\slsf
m}-normality (\ie $\Gamma$-normality in our sense) was given by
Fujimoto in \cite{Fu}.

\newprg[INT.rat-con]{Rational connectivity of the exceptional
components of the limit} Strong convergence obviously implies the
weak one and the latter implies the gamma-convergence, see Remark
\ref{weak-gamma}:
\begin{equation}
\eqqno(implic) {\sf s}\text{-convergence} \Longrightarrow {\sf
w}\text{-convergence} \Longrightarrow \Gamma\text{- convergence}.
\end{equation}
Our second principal task in this paper is to understand what
obstructs a weakly/gamma converging sequence to converge strongly.
The problem is that by Theorem \ref{str-cycl} the volumes of graphs
of a strongly converging sequence are uniformly bounded over
compacts in the source. When dimension $n$ of the source $U$ is two
and $X$ is K\"ahler the volumes of the graphs of a weakly converging
sequence are still bounded, see \cite{Iv2}. The same is true if is
$X$ an arbitrary compact complex surface (and again $\dim U = 2$),
see \cite{Ne}. We shall say more about this in section \ref{VOL}.
But this turns out {\slsf not to be the case} starting from
dimension three, \ie the volumes of graphs of a weakly converging
sequence can diverge to infinity over compacts of $U$. Via
\eqqref(monge1) this turns out to be a geometric counterpart of a
well known discontinuity of Monge-Amp\`ere masses, see Example
\ref{rash} in section \ref{PN}. Nevertheless  for a sequence
$\Gamma_{f_k}$ of $\Gamma$-converging meromorphic graphs we can
consider the Hausdorff limit $\hat\Gamma$ (its always exists after
taking a subsequence). Set $\Gamma \deff \overline{\hat\Gamma
\setminus \Gamma_f}$, where $\Gamma_f$ is the graph of the limit map
$f$, and call $\Gamma$ {\slsf a bubble}. For $a\in \gamma\deff
\pr_1(\Gamma)$ set $\Gamma_a \deff \pr_2(\pr_1^{-1}(a)\cap \Gamma)$,
here $\pr_1$ and $\pr_2$ are natural projections, see section
\ref{MER-MAP}. We prove the following statement.

\begin{nnthm}
\label{rat-con0} Let $X$ be a  disk-convex Gauduchon manifold and
let $f_k:U\to X$ be a weakly converging sequence of meromorphic
mappings which doesn't converge strongly. Then for every point $a\in 
\gamma$ the fiber $\Gamma_a$ is
rationally connected. If $X$ is, moreover, projective then the same
is true also for a $\Gamma$-converging sequences.
\end{nnthm}
Here by saying that a closed subset $\Gamma_a$ of a complex manifold
is {\slsf rationally connected} we mean that every two distinct
points $p,q\in \Gamma_a$ can be connected by a chain of rational
curves which is entirely contained in $\Gamma_a$, see section
\ref{RAT-CON} for more details.

\newprg[INT.fatou]{Fatou sets of meromorphic self-maps}
Families of a special interest are the families of {\slsf iterates}
$f^n \deff f\circ ... \circ f$ of some fixed meromorphic self-map of
a {\slsf compact} complex manifold $X$. The maximal open subset $X$
where $\{f^n\}$ is relatively compact is called the Fatou set of
$f$. Depending on the sense of convergence that one wishes to
consider one gets different Fatou sets: strong, weak or gamma Fatou
sets. We denote them as $\Phi_s$, $\Phi_w$ and $\Phi_{\Gamma}$ respectively,
their dependance on $f$ will be clear from the context.

\begin{nncorol}
\label{fatou} Let $f$ be a meromorphic self-map of a compact complex
surface. Then the weak Fatou set $\Phi_w$ of $f$ is locally
pseudoconvex. If $\Phi_s$ is different from $\Phi_w$ then:

\smallskip ${\sf a)}$ $X$ is bimeromorphic to $\pp^2$;

\smallskip ${\sf b)}$ $\Phi_w = X\setminus C$, where $C$ is a
rational curve in $X$;

\smallskip ${\sf c)}$ the weak limit of any weakly converging
subsequence $\{f^{n_k}\}$ of iterates is a degenerate

\quad map of $X$ onto $C$.
\end{nncorol}

It should be pointed out that our Fatou sets are different from the
Fatou sets as they were considered in \cite{FS}. In \cite{FS} the
Fatou set of $f$ is the maximal open subset $\Phi$ of
$X\setminus\overline{\bigcup_{n=0}^{\infty}f^{-n}(I_f)}$ where the
family $\{f^n\}$ is equicontinuous (remark that on
$X\setminus\overline{\bigcup_{n=0}^{\infty}f^{-n}(I_f)}$ all
iterates are holomorphic). If, for example, $f:\pp^2\to \pp^2$ is
the Cremona transformation $[z_0:z_1:z_2]\to [z_1z_2:z_0z_2:z_0z_1]$
then $\Phi_s = \Phi_w = \Phi_{\Gamma} = \pp^2$ but $\Phi =
\pp^2\setminus \{ \text{ three lines } \}$. In subsection
\ref{FATOU.exmp} an example of higher degree and with an interesting
dynamics on the indeterminacy set is given. This is one more instance
which shows how crucially can
change a picture when the notion of convergence changes.

\medskip\noindent{\slsf Notes. 1.} Let us make a final note about the goals
of this paper. On our opinion the most interesting information about
a converging sequence of meromorphic mappings is concentrated near
the ``limit`` of their indeterminacy sets.  We describe the most
reasonable (in our opinion) notions of convergence of meromorphic
mappings and conclude that the {\slsf weak} one is the most
appropriate. At the same time we detect that if a {\slsf
weakly/gamma} converging sequence doesn't converge {\slsf strongly}
then this imposes very serious restrictions on the target manifold
(it is forced to contain many rational curves). In some cases (ex.
iterations) this puts strong constraints also on the sequence
itself.

\smallskip\noindent{\slsf 2.} Domains of convergence of holomorphic
functions of several variables were, probably, for the first time
considered by G. Julia in \cite{J}. In \cite{J} and then in
\cite{CT} it was proved that these domains are (in some sense)
pseudoconvex. Domains of convergence of meromorphic functions of
several variables were studied in \cite{Sa} and then in \cite{Ru}.
In these early papers convergence was understood as holomorphic (\ie
uniform) convergence outside of the union of indeterminacy sets of
meromorphic mappings in question.

\medskip\noindent{\slsf Acknowledgement.} We are  grateful
to Alexander Rashkovskii for explaining to us the Example \ref{rash} with
unbounded Monge-Amp\`ere masses.

\newsect[MER-MAP]{Topologies on the space of meromorphic mappings}

\newprg[MER-MAP.man-map]{Complex manifolds and meromorphic mappings}

Our manifolds will be Hausdorff and countable at infinity if the
opposite is not explicitly stated. We shall also everywhere suppose
that they are {\slsf disk-convex}.

\begin{defi}
\label{disk-conv} A complex manifold $X$ is called disk-convex if
for every compact $K\comp X$ there exists a compact $\hat K$ such
that for every $h\in \calo (\Delta ,X)\cap \calc(\bar\Delta , X)$
such that $h(\d\Delta)\subset K$ one has $h(\bar\Delta)\subset \hat
K$.
\end{defi}

\noindent The minimal such $\hat K$ is called the {\slsf disk envelope}
of $K$. Let $X$ be
equipped with some Hermitian metric $h$. By $\omega_h$ denote the
$(1,1)$-form canonically associated with $h$. We say that the metric
$h$ is $d$-{\slsf closed} or {\slsf K\"ahler} if $d\omega_h=0$. We
say that $h$ is {\slsf pluriclosed} or {\slsf
Gauduchon} if $dd^c\omega_h = 0$. In \cite{Ga} it was proved
that on a compact complex surface every Hermitian metric
is conformally equivalent to the unique $dd^c$-closed one.
\begin{rema}\rm
\label{gaud}
We shall need only the existence of such metric forms on
compact complex surfaces and this can be proved by duality: non
existence of a positive $dd^c$-closed $(1,1)$-form is equivalent to
the existence of a non-constant plurisubharmonic function. The latter
on a compact complex manifold is impossible.
\end{rema}
We also fix some metric form $\omega_1$ on $U$. In the case of a
polydisk $U = \Delta^n$ we will work with the  standard Euclidean
metric $e$. The associated form will be denoted by $\omega_e = dd^c
\Vert z\Vert ^2 = \frac{i}{2}\sum_{j=1}^ndz_j\wedge d\bar z_j $. By
$\pr_1: U \times X \longrightarrow U$ and $\pr_2: U\times
X\longrightarrow X$ denote the projections onto the first and second
factors. On the product $U\times X$ we consider the metric form
$\omega = \pr_1^{\ast }\omega_1 + \pr_2^{\ast }\omega_h$.

\smallskip A meromorphic mapping $f$ between complex manifolds $U$
and $X$ is defined by an {\slsf irreducible}
analytic subset $\Gamma_f\subset U\times X$ such that

\begin{itemize}
\item {\it the restriction $\pr_1|_{\Gamma_f} : \Gamma_f \to U$ of
the natural projection to $\Gamma_f$ is a proper
modification, \ie is proper and generically one to one.}
\end{itemize}

$\Gamma_f$ is called the {\slsf graph} of $f$. Due to the
irreducibility of $\Gamma_f$ and the Remmert proper mapping theorem
the set of points over which $\pr_1$ is not one to one is an
analytic subset of $U$ of codimension at least two. This set is
called the set of points of indeterminacy of $f$ and is usually
denoted as $I_f$. Therefore an another way to define a meromorphic
mapping $f$ between complex manifolds $U$ and  $X$ is by considering
a holomorphic map $f:U\setminus A\to X$, where $A$ is an analytic
subset of $U$ of codimension at least two, such that the closure
$\Gamma_f$ of its graph is an analytic subset of the product
$U\times X$ satisfying the  condition above. Remark that the
analyticity of the closure of the holomorphic graph is not
automatic. Think about the natural projection $f:\cc^2\setminus \{0\}\to
\cc^2\setminus \{0\}/z\sim 2z$ of $\cc^2\setminus \{0\}$ onto a Hopf
surface. The properness of the restriction of the projection $\pr_1$
to the closure is, unless $X$ is disk convex, not automatic too.

\smallskip
The volume of the graph $\Gamma _{f}$ of a meromorphic mapping $f$ is given by

\begin{equation}
\eqqno(gr-vol1)
n!\vol (\Gamma_f) = \int\limits_{\Gamma_f}\omega^n =
\int\limits_{\Gamma_f}(\pr_1^{\ast }\omega_1 + \pr_2^{\ast
}\omega_h)^n =\int\limits_U\big(\omega_1 + f^{\ast }\omega_h\big)^n,
\end{equation}
where $n=\dim U$.

\begin{rema} \rm
\label{mer-rem} Let us make a few remarks concerning the notion of a
meromorphic mapping.

\smallskip\noindent{\sl a)} If $V$ is a subvariety of $U$  such that
$V\not\subset I_f$ then the restriction $f|_V$ of $f$ to $V$ is
defined by taking as its graph $\Gamma_{f|_V}$ the {\slsf
irreducible component} of the intersection $\Gamma_f\cap (V\times
X)$ which projects onto $V$ generically one to one. Therefore
$\Gamma_{f|_V}\subset \Gamma_f\cap (V\times X)$ and the inclusion
here is proper in general. The {\slsf full image} of a set $L\subset U$ under
$f$ is defined as $f[L]\deff \pr_2\left(\Gamma_f\cap [L\times X]\right)$.

\smallskip\noindent{\sl b)} It is probably worth to notice that $x\in I_f$
if and only if $\dim f[x]\ge 1$. This follows from the obvious
observation that $I_f = \pr_1\left(\{(x_1,x_2)\in \Gamma_f:
\dim_{(x_1,x_2)} \pr_1|_{\Gamma_f}^{-1}(x_1)\ge 1\}\right)$.

\smallskip\noindent{\sl c)} If $\dim V = 1$ then the irreducible component
of $\Gamma_f\cap (V\times X)$
which projects onto $V$ is a curve. Since the projection is
generically one to one it is on to one everywhere and therefore the
restriction $f|_V$ is necessarily {\slsf holomorphic}.

\smallskip\noindent{\sl d)} Let us give the sense to $f^{\ast }\omega_h$
in the formula \eqqref(gr-vol1). The first integral there has
perfectly sense since we are integrating a smooth form over a
complex variety. Denote by $I_f^{\eps}$ the $\eps$-neighborhood of
the indeterminacy set $I_f$ of $f$. Then \eqqref(gr-vol1) shows that
the limit
\begin{equation}
\eqqno(gr-vol2)
\lim_{\eps\to 0} \int\limits_{U\setminus\bar
I_f^{\eps}}\big(\omega_1 + f^{\ast }\omega_h\big)^n = \lim_{\eps\to
0} \int\limits_{U\setminus\bar
I_f^{\eps}}\sum_{p=0}^nC_n^p\omega_1^{n-p}\wedge f^{\ast
}\omega_h^{p}
\end{equation}
exists. Therefore all $f^{\ast }\omega_h^{p}$ are well defined on
$U$ as positive currents.
\end{rema}

\newprg[MER-MAP.an-cycle]{Analytic cycles and currents}

Before turning to the notions of convergence of meromorphic mappings
let us recall the natural topologies on the space of analytic
subsets of a complex manifold.

\smallskip Recall that an analytic cycle of dimension $r$ in a complex
manifold $Y$ is a formal sum $Z=\sum_jn_jZ_j$, where $\{ Z_j\} $ is
a locally finite sequence of reduced analytic subsets  of pure
dimension $r$ and $n_j$ are positive integers called multiplicities
of $Z_j$. The set $\vert Z\vert :=\bigcup_jZ_j$ is called the {\slsf
support} of $Z$. In our applications $Y$ will be $U\times X$ and $r$
will be the dimension $n=\dim U$. By a coordinate chart adapted to $Z$ we
shall understand a relatively compact open set $V$ in $Y$ such that
$V\cap \vert Z\vert \not=\emptyset $ together with a biholomorphism
$j$ of $V$ onto a neighborhood $\tilde V$ of $\bar\Delta^r\times
\bar\Delta^q$ in $\cc^{r+q}$,  $r+q=\dim Y$, such that
$j^{-1}(\bar\Delta^r \times \partial \Delta^q)\cap \vert Z\vert =
\emptyset $. We shall denote such chart by $(V,j)$. The image
$j(Z\cap V)$ of the cycle $Z\cap V$ under biholomorphism $j$ is the
image of the underlying analytic set together with multiplicities.
Following Barlet and Fujiki, see \cite{Ba} and \cite{Fj}, we call
the quadruple $E= (V,j,\Delta^r, \Delta^q)$ a  {\slsf scale} adapted
to $Z$.

\smallskip If $\pr :\cc^r \times \cc^q\to \cc^r $ is the natural projection,
then the restriction $\pr\mid_{j(Z\cap V)}:j(Z\cap V)\to \Delta^r $
is a branched covering of degree say $d$. This branched covering
defines in a natural way a holomorphic mapping
$\phi_{j,Z}:\Delta^r\to \sym^d (\Delta^q)$ to the $d$-th symmetric
power of $\Delta^q$ by setting
$\phi_{j,Z}(z')=\left\{(\pr\mid_{j(Z\cap V)} )^{-1}(z')\right\}$.
The latter denotes the
unordered set of all preimages of $z^{'}$ under the projection in
question. This construction, due to Barlet, allows to represent a
cycle $Z\subset Y$ by a set of holomorphic maps
$\phi_{j_{\alpha},Z}:\Delta^r\to \sym^d (\Delta^q)$, where
$\{(V_{\alpha},j_{\alpha})\}$ is some open covering of $|Z|$ by
adapted coordinate charts.

\begin{defi}
One says that $Z_k$ converges to $Z$ in the topology of cycles if
for every coordinate chart $(V,j)$ adapted to $Z$ there exists $k_0$
such that $\forall k\ge k_0$ this chart will be adapted to $Z_k$ and
the sequence of corresponding holomorphic mappings $\phi_{j,Z_k}$
converge to $\phi_{j,Z}$ uniformly on $\Delta^r$.
\end{defi}

This defines a metrizable topology on the space $\calc_r(Y)$ of
$r$-cycles in $Y$. This topology is equivalent to the {\slsf topology of
currents}: $Z_k\to Z$ if for any continuous $(r,r)$-form $\chi$ with
compact support one has
\[
 \int\limits_{Z_k}\chi \to \int\limits_{Z}\chi ,
\]
see \cite{Fj}. It is also equivalent to the {\slsf Hausdorff} topology under an
additional condition of {\slsf boundedness of volumes}. Recall that
the Hausdorff distance between two subsets $A$ and $B$ of a metric
space $(Y,\rho )$ is a number $\rho (A,B) = \inf \{ \varepsilon :
A^{\varepsilon } \supset B, B^{\varepsilon }\supset A\} $. Here by
$A^{\varepsilon }$ we denote the $\varepsilon $-neighborhood of the
set $A$, i.e. $A^{\varepsilon } = \{y\in Y: \rho (y,A) < \varepsilon
\}$.

\smallskip
Now, $Z_k\to Z$ if and only if for every compact $K\comp Y$ there
exists $C_K>0$ such that $\vol_{2r}(Z_k\cap K)\le C_K$  and $Z_k\cap
K \to Z\cap K$ with respect to the Hausdorff distance. This
statement is the content of the Harvey-Shiffman's generalization of
Bishop's compactness theorem. For the proof see \cite{HS}. We denote
the space of $r$-cycles on $Y$ endowed with the topology described
as above by $\calc^{\loc}_r(Y)$.

\newprg[MER-MAP.s-conv]{Strong convergence of meromorphic mappings}

Let $\{ f_k\} $ be a sequence of meromorphic mappings of a complex
manifold $U$ to a complex manifold $X$.

\begin{defi}
\label{str-def} We say that $f_k$ converge {\slsf strongly} to a
meromorphic map $f:U \to X$ ({\slsf s}-converge) if the sequence of
graphs $\Gamma_{f_k}$ converge over compacts to $\Gamma_f$ in
Hausdorff metric, \ie for every compact $K\comp U$ one has
$\Gamma_{f_k}\cap (K\times X)\xrightarrow{H} \Gamma_{f}\cap (K\times
X)$.
\end{defi}

Now let us prove Theorem \ref{str-cycl} from the Introduction, \ie that
Hausdorff convergence in the case of graphs implies the boundedness
of volumes (over compacts) and therefore the convergence in the
{\slsf topology of cycles}. Let us underline at this point that in
this theorem one doesn't need to suppose anything on the target
manifold $X$.

\medskip\noindent{\slsf Proof of Theorem \ref{str-cycl}.}
The reason why Hausdorff convergence of graphs implies their
stronger convergence in the topology of cycles is that, being the
graphs, the analytic cycles $\Gamma_{f_k}$ converge to $\Gamma_f$
with {\slsf multiplicity} one. Now let us give the details. Let $a
\in U\setminus I_f$ be a regular point of $f$ and set $b=f(a)$. Then
we can find neighborhoods $D_1\ni a$ and $D_2\ni b$ biholomorphic to
$\Delta^n$, $n=\dim U$ and $\Delta^p$, $p=\dim X$ respectively such
that $\Gamma_f\cap \left(\bar D_1 \times \d D_2\right)=\emptyset$.
In particular $V=D_1\times D_2$ is an adapted chart for $\Gamma_f$,
let $(V, j, \Delta^n, \Delta^p)$ be a corresponding scale. Here
$j:V\to \Delta^n\times \Delta^p$ is some biholomorphism. From Hausdorff
convergence of $\Gamma_{f_k}$ to $\Gamma_f$ we see that for $k\gg1$
$\Gamma_{f_k}\cap  \left(\bar D_1 \times \d D_2\right)=\emptyset$.
Therefore $\Gamma_{f_k}\cap  \left(D_1 \times D_2\right) \to D_1$ is
a ramified covering (\ie is proper) of some degree $d_k$. But
$\Gamma_{f_k}$ is one to one over a generic point of $D_1$.
Therefore $d_k=1$ and $\Gamma_{f_k}\cap V$ converge to $\Gamma_f\cap
 V$ as graphs (in particular as cycles). We proved that $f_k$
converge to $f$ on compacts of $U\setminus I_f$ as holomorphic
mappings.

\smallskip Let now $a\in I_f$ and take some $b\in f[a]$. As above
take a neighborhood $V = D_1\times D_2\cong \Delta^{n}\times
\Delta^p$ of $(a,b)$, where $a=0$ and $b=0$ in these coordinates.
Denote by $(w^{'},w^{''})$ the coordinates in $\Delta^{n}\times
\Delta^p$. Perturbing the slope of coordinate $w^{''}$ we can
suppose that $(\{0\}\times \Delta^p)\cap \Gamma_f$ has $0$ as its
isolated point.

\begin{rema} \rm
After perturbation of the slope of $w^{''}$ the decomposition $j(V)
= \Delta^n\times \Delta^p$ will not correspond to the decomposition
$U\times X$.
\end{rema}

For sufficiently small $\eps > 0$ we take  polydisks
$\Delta^{n}_{\eps}$ and $\Delta^p_{\eps}$  in (perturbed)
coordinates (actually only $w^{''}$ needs to be perturbed). We get
an adapted chart for $\Gamma_f$ which possed the following property:

\medskip\sli $\tilde V \deff j(V)$ writes as $\tilde V =
\Delta^{n}\times \Delta^p$, $p = \dim X$.

\smallskip\slii Local coordinates $(w^{'},w^{''})$ of
$\Delta^{n}\times \Delta^p$ enjoy the property that $z^{'}\deff w^{'} \in
\Delta^{n}$ is a local

\quad   coordinate in $U$ (but $w^{''}$ is not a local coordinate on $X$).

\smallskip\sliii $j(\Gamma_{f}\cap V) \to \Delta^n$ is a ramified covering
of degree $d\ge 1$.

\medskip Again from Hausdorff convergence of $\Gamma_{f_k}$ to $\Gamma_f$
we get that for all $k\gg1$ the intersection $j(\Gamma_{f_k}\cap V )$
is a ramified covering of $\Delta^n$ of degree $d_k$. Obviously
$d_k\ge d$ for $k\gg1$. If for some subsequence $d_k>d$ we shall get
a contradiction as follows. In that case some irreducible component
of $\Gamma_f\cap V$ will be approached by $\Gamma_{f_k}\cap V$ at
least doubly. Let $\Gamma$ stands for this irreducible component.
Since $\dim \left[\Gamma_f\cap (I_f\times X)\right]\le n-1$ (by
irreducibility of $\Gamma_f$) we see that $\Gamma_{f_k}$ multiply
approach every compact of $\Gamma\setminus (I_f\times X)$. Take a
point $c\in \Gamma\setminus (I_f\times X)$ having a relatively
compact neighborhood $W\subset \Gamma\setminus (I_f\times X)$ such
that $\pr_1|_W:W\to W_0$ is biholomorphic, \ie $W$ is the graph of
$f$ over $W_0\comp U\setminus I_f$. Now it is clear that
$\Gamma_{f_k}\cap (W_0\times X)$ cannot approach $\Gamma_f\cap
(W_0\times X) = W$ with multiplicity more than one, because
$\Gamma_{f_k}$ is a graph of a holomorphic map over $W_0$ for
$k\gg1$.

\smallskip We proved that every irreducible branch of $\Gamma_f\cap V$
the graphs $\Gamma_{f_k}$ approach with multiplicity one. Therefore
$j(\Gamma_{f_k}\cap V) \to \Delta^n$ is a ramified covering of the
same degree $d$ as $j(\Gamma_{f}\cap V) \to \Delta^n$ for $k\gg 1$.
This proves at a time that $\Gamma_{f_k}$ converge to $\Gamma_f$ in
the topology of cycles and that their volumes are uniformly bounded.

\smallskip\qed

\smallskip
Strong convergence has some nice features, one was mentioned in the
Introduction. Moreover, as it is explained in \cite{Iv4}, strong topology
is natural in studying fix points of meromorphic self-mappings of
compact complex manifolds. But domains of strong convergence and
strong normality are quite arbitrary.
We shall explain this in more details. Let $\calf$ be a family of meromorphic
mappings from a complex manifold $U$ to a disk convex complex manifold $X$.

\begin{defi}
\label{norm-set}
The {\slsf set of normality} of $\calf$ is the maximal open subset
$\caln_{\calf}$ of $U$ such that $\calf$ is relatively compact on
$\caln_{\calf}$. If $\calf = \{f_k\}$ is a sequence then the {\slsf
set of convergence} of $\calf$ is the maximal open subset of $U$
such that $f_k$ converge on compacts of this subset.
\end{defi}

To be relatively compact in this definition means that from every
sequence of elements of $\calf$ one can extract a converging on
compacts subsequence. The sense of convergence (strong, weak or
other) should be each time specified.

\begin{exmp} \rm
\label{norm-exmp}
{\bf 1.} Let $X$ be a Hopf three-fold $X\deff \cc^3\setminus \{0\}/z\sim 2z$.
Denote by $\pi:\cc^3\setminus \{0\}\to X$ the canonical projection.
Let $D\comp \cc^2$ be any bounded domain. Take a sequence $\{a_n\}\subset D$
accumulating to every point on $\d D$. Let $g_n:\cc^2\to \cc^3$ be defined
as $g_n(z) = (z-a_n, 1/n)$. Set $f_n\deff \pi\circ g_n$. Then the set
of normality of $\calf = \{f_n\}$ has $D$  as one of its connected
components.

\begin{rema}\rm
For an analogous example with $X$ projective see  Example 4 from \cite{Iv2}.
\end{rema}

\noindent{\bf 2.} The same example is instructive when understanding the
notion of weak convergence. Take a converging to zero sequence $a_n$.
Then $f_n$ from this example will converge on compacts of $\cc^2\setminus \{0\}$
but the limit will not extend to zero meromorphically. I.e., $f_n$ will not
converge weakly in any neighborhood of the origin.
\end{exmp}

\newsect[PSCN]{Pseudoconvexity of sets of normality}

\newprg[PSCN.w-conv]{Weak convergence and  proof of Theorem \ref{propa}}

In view of such examples a weaker notion of convergence for meromorphic
mappings was introduced in \cite{Iv2}. Let $f\in \calm (U,X)$ be a
meromorphic map from $U$ to $X$ and let $\{f_k\}\subset \calm (U,X)$
be a sequence of meromorphic mappings.

\begin{defi}
\label{weak-def}
We say that $f_k$ converge {\slsf weakly} to $f$ ({\slsf
w}-converge) if there exists an analytic subset $A$ in $U$ of
codimension at least two such that $f_k$ converge strongly to $f$ on
$U\setminus A$.
\end{defi}

\begin{rema} \rm
\label{weak-rem} $f_k$ converge weakly to $f$ if and only if for
every compact of $U\setminus I_f$ all $f_k$ are holomorphic in a
neighborhood of this compact for $k$ big enough and uniformly
converge there to $f$ as holomorphic mappings. Indeed, let $A$ be
the minimal analytic set of codimension $\ge 2$ such that $f_k$
converge strongly to $f$ on $U\setminus A$. Then $A$ must be
contained in $I_f$ because if there exists a point $a\in A\setminus
I_f$ then $f$ is holomorphic in some neighborhood $V\ni a$ and then,
by Rouch\'e Principle of \cite{Iv2} $f_k$ for $k\gg 1$ are holomorphic
on compacts in $V\setminus A$ and converge uniformly (on compacts)
to $f$ there. From here and the fact that $\codim A\ge 2$ one easily
gets that $f_k$ are holomorphic on compacts in $V$ and converge to
$f$.
\end{rema}

Now let us turn to the sets of weak convergence/normality. Sets of strong
normality obviously are well defined, \ie they do exist. The existence of sets
of weak normality was proved in \cite{Iv2}, see Corollary 1.2.1a.

\begin{rema} \rm
In the formulation of this Corollary the Author of \cite{Iv2} speaks about
``weak convergence'' but the proof is about ``weak normality``.
\end{rema}

\smallskip Domains of weak convergence of meromorphic mappings turn to be
pseudoconvex for a large class of target manifolds. This follows
from the "mutual propagation principle" stated in Theorem \ref{propa}
in the Introduction. Let us give a proof of it.

\medskip\noindent{\slsf Proof of Theorem \ref{propa}.}
Let us prove the part (b) first.

\medskip\noindent{\slsf Step 1. Extension of the limit.} First of all by the main
result of \cite{Iv3} every
meromorphic map $f:U\to X$ extends to a meromorphic map $f:\hat U \setminus
A\to X$, where $A$ is closed, complete $(n-2)$-polar subset of $\hat U$ of
Hausdorff $(2n-3)$-measure zero.
In more details that means that for every point $a\in A$ there exists a coordinate
neighborhood $V\cong \Delta^{n-2}\times \bb^2$ of $a=0$ such that
$A\cap (\Delta^{n-2}\times \d \bb^2) = \emptyset$ and for every
$z^{'}\in \Delta^{n-2}$ the intersection $A_{z^{'}}\deff A\cap \bb^2_{z^{'}}$
is a zero dimensional complete pluripolar subset of $\bb^2_{z^{'}}
\deff \{z^{'}\}\times \bb^2$. Here $\bb^2$ stands for the unit ball in $\cc^2$.
Moreover, if $A\not=\emptyset $ then $f(\sss^3_{z^{'}})$ is not homologous to
zero in $X$. Here $\sss^3_{z^{'}}=\d \bb_{z^{'}}$ is the standard three-dimensional
sphere in $\cc^2$.

\smallskip Let $U'$ be the maximal open subset of $\hat U\setminus (I_f\cup A)$ such that
$f_k$ converge to $f$ on compacts on $U'$ as holomorphic mappings.

\medskip\noindent{\slsf Step 2. $U'$ is locally pseudoconvex in $\hat U\setminus  (I_f\cup A)$.}
If not  then by Docquier-Grauert criterion, see \cite{DG},
there would exist a point $b\in \d U'\setminus  (I_f\cup A)$ and
a Hartogs figure $h:H^n_{\eps}\to U^{'}$ imbedded to $U'$  such
that the image $h(\Delta^n)$ of the corresponding polydisk contains
$b$. All this is local and therefore we can assume that
$h(\Delta^n)$ is relatively compact in $U\setminus  (I_f\cup A)$.
Pulling back $f_k$ and $f$ to $\Delta^n$ we arrive to contradiction
as follows. By the Theorem of Siu, see \cite{Si2}, there exists a
Stein neighborhood $V$ of the graph of $f\circ h$
in $\Delta^n\times X$. Since for every compact $K\comp H^n_{\eps}$
we have that the graph of $f_k|_K$ is contained in $V$ we conclude
the same for every compact of $\Delta^n$. Now $f_k\circ h$ converge
to $f\circ h$ on compacts in $\Delta^n$ as holomorphic mappings. But
that mean that they converge also around the preimage of $b$.
Contradiction. Since $\hat U$ was supposed to be the envelope of
holomorphy of $U$ we obtain that $U' = \hat U\setminus  (I_f\cup A)$.

\smallskip\noindent{\slsf Step 3. Removing $A$.} Suppose $A$ is non-empty. Take a sphere
$\sss^3_{z^{'}}$ as described in {\slsf Step 1} for some point $a\in A$. Using
the fact that $I_f$ is of codimension $\ge 2$ we can take  $\sss^3_{z'}$ not
to intersect $I_f$ as well. I.e., $\sss^3_{z'}\subset U'$.  $f_k(\sss^3_{z'})$ is homologous to
zero in $X$, because $f_k$ meromorphically extends to the corresponding $\bb^2_{z'}$.
Moreover $f_k$ converge to $f$ in a neighborhood of $\sss^3_{z'}$. This implies
that $f(\sss^3_{z'})$ is also homologous to zero and therefore $A$ should be empty.
Contradiction. Part (b) is proved.

\medskip The proof of (a) is a particular case of the {\slsf Step 2} of the proof
of part (b).

\smallskip\qed

\begin{rema} \rm
We gave a proof of Theorem \ref{propa} here because the proof of
an analogous statement in \cite{Iv2} uses a stronger extension claim
from the subsequent paper \cite{Iv3}. Namely the Author claimed that
$A$ appearing in the {\slsf Step 1} of the proof is {\slsf analytic}
of codimension two. This was not achieved in \cite{Iv3} (and is not
clear for us up to know). Therefore we find necessary to remark that
vanishing of $(2n-3)$-dimensional measure of $A$ together with homological
characterization of the obstructions for the meromorphic extension is,
in fact, sufficient for our particular task here.
\end{rema}

\newprg[PSCN.gamma]{Gamma convergence of meromorphic mappings}

Let again $f_k$ be a sequence of meromorphic mappings between complex
manifolds $U$ and $X$, the last is supposed to be disk-convex. Let
$f\in \calm (U,X)$ be a meromorphic map.

\begin{defi}
\label{gamma-def}
We say that $f_k$ $\Gamma$-converge to  $f$ if:

\smallskip\sli there exists an analytic subset $A\subset U$ such that
$f_k$ strongly converge to $f$ on $U\setminus A$;

\smallskip\slii for every divisor $H$ in $X$, such that $f(U)\not\subset H$
and every compact $K\comp U$ the

volumes of $f_k^*H\cap K$ counted with multiplicities are
uniformly bounded for $k\gg 1$.
\end{defi}

\begin{rema} \rm
\label{weak-gamma} This notion is strictly weaker than the weak
convergence because $A$ can have components of codimension one, and
remark that the item (\slii is automatically satisfied by a weakly
converging sequence, because divisors $f^*H$ extend from $U\setminus
A$ to $U$ and if they have bounded volume on compacts of $U\setminus
A$ then the same is true on compacts of $U$. All this obviously
follows from the ingredients involved in the proof of Bishop's
compactness theorem, see \cite{Bi} or \cite{St}. It might be
convenient to add to $A$ the indeterminacy set of $f$ and then, see
Remark \ref{weak-rem},  $f_k$ will converge to $f$ uniformly on
compacts of $U\setminus A$ as holomorphic mappings.

\end{rema}

\begin{exmp} \rm
\label{exp}
{\slsf a)} Consider the following sequence of holomorphic mappings $f_k:\Delta \to
\pp^1$:
\begin{equation}
\eqqno(exp)
f_k : z \to \left[1: 1 + \frac{1}{z} +  ... + \frac{1}{z^kk!}\right] =
\left[z^k : z^k + z^{k-1} + ... + \frac{1}{k!}\right].
\end{equation}
It is clear that $f_k$ converges on compacts of $\Delta\setminus \{0\}$ to
$f(z) = [1:e^{\frac{1}{z}}]$ but, as it is clear from the second expression
in \eqqref(exp) the preimage counting with multiplicities of the divisor
$H= \{Z_0=0\}$ is $k[0]$ (here $[Z_0:Z_1]$ are homogeneous coordinates in
$\pp^1$), \ie has unbounded volume. And indeed, this
sequence should not be considered as converging one, because its limit is
not holomorphic on $\Delta$.

\smallskip\noindent{\slsf b)} Set $f_k(z) = [z: z-\frac{1}{k}]:\Delta\to \pp^1$. This
sequence clearly converges to the constant map $f(z) = [z:z] = [1:1]$ on compacts of $\Delta
\setminus \{0\}$. Moreover, the preimage of any divisor $H = \{ P(z_0,z_1)=0\}$ under
$f_k$ is $\{ z\in \Delta : P(z,z-\frac{1}{k})=0\}$, \ie is a  set of points,
uniformly bounded in number counting with multiplicities.
Therefore this sequence $\Gamma$-converge (but doesn't converge weakly).
\end{exmp}

\begin{exmp} \rm
\label{rutish}
Consider the following sequence of meromorphic functions on $\Delta^2$ (\ie meromorphic
mappings to $\pp^1$):
\[
f_k(z_1,z_2) = [z_1:2^{-k}z_2^k].
\]
The limit map is constant $f(z) = [1:0]$. $f_k$ converge to $f$
strongly (uniformly in fact) on compacts of $\Delta^2\setminus
\{z_1=0\}$. If $[Z_0:Z_1]$ are homogeneous coordinates in $\pp^1$
then the preimage of the divisor $[Z_1=0]$ is $k[z_2=0]$, \ie this
sequence doesn't converge even in $\Gamma$-sense on $\Delta^2$.
\end{exmp}

\begin{rema} \rm
Examples \ref{exp} (a) and \ref{rutish} are examples of converging
outside of an analytic set of codimension one sequences which are
not $\Gamma$-converging. In the first case the limit doesn't extend
to the whole source, in the second it does. Convergence of
meromorphic mappings of this type was introduced and studied by
Rutishauser in \cite{Ru}.
\end{rema}

If in Definition \ref{norm-set} the underlying convergence is
$\Gamma$-convergence we get the corresponding notions of a
convergence/normality set. Let us conclude this general discussion
with the following

\begin{prop}
\label{gamma-norm} Let $X$ be a disk-convex Gauduchon manifold. Then
the sets of $\Gamma$-convergence/normality of meromorphic mappings
with values in $X$ are locally pseudoconvex.
\end{prop}
\proof We shall prove the statement for the sets of
$\Gamma$-normality, the case of sets of convergence obviously
follows. Let $D$ be the maximal open subset of $U$ where the family
$\calf$ is $\Gamma$-normal. Suppose that $D$ is not pseudoconvex.
Then by Docquier-Grauert criterion, see \cite{DG}, there exists an
imbedding $h:H^n_{\eps} \to D$ of a Hartogs figure into $D$ such
that $h$ extends to an immersion of the polydisk to $U$ with
$h(\Delta^n)\cap \d D\not= \emptyset$. Recall that Hartogs
figure is the following domain
\begin{equation}
\eqqno(hart-fig)
H^n_{\eps} \deff \left(\Delta^{n-1}_{\eps} \times \Delta\right)\cup \left(\Delta^{n-1}\times
A_{1-\eps , 1}\right),
\end{equation}
where $A_{1-\eps , 1} \deff \Delta\setminus \bar\Delta_{1-\eps}$ is an annulus.
Let us pull-back our family
to $\Delta^n$ by $h$ and therefore without loss of generality we can
suppose that $U=\Delta^n$, $\calf$ is a family of meromorphic
mappings from $\Delta^n$ to $X$, $H^n_{\eps}\subset D\subset
\Delta^n$ is the set of $\Gamma$-normality of $\calf$ such that
$D\not= \Delta^n$.

\smallskip That means that there exists a sequence $\{f_k\}\subset \calf$,
which converges on $D$ but doesn't not $\Gamma$-converge on compacts in $\Delta^n$. Let us see
that this is impossible. Let $f:D\to X$ be the $\Gamma$-limit of
$f_k$. Denote by $A$ the analytic set in $D$  such that $f_k$
converge to $f$ on compacts of $D\setminus A$. By \cite{Iv3} $f$
extends to $\Delta^n\setminus S$, where $S$ is closed
$(n-2)$-complete polar subset of $\Delta^n$. Let $A'$ be the pure
$(n-1)$-dimensional part of  $A$. By the theorem of Grauert we have
two cases.

\smallskip\noindent{\slsf Case 1. The envelope of holomorphy of $D\setminus A'$ is $\Delta^n$.}
In that case the Theorem \ref{propa} is applicable with $U=D\setminus A'$ and
$\hat U =\Delta^n$ and gives us the weak (and therefore $\Gamma$) convergence of $f_k$
on $\Delta^n$.

\smallskip\noindent{\slsf Case 2. $A'$  extends to a hypersurface $\tilde A$ in $\Delta^n$ and
$\Delta^n\setminus \tilde A$ is the envelope of holomorphy of
$D\setminus A'$.} In that case again by Theorem \ref{propa}
$f_k$ weakly converge to $f$ on $\Delta^n\setminus \tilde A$.
$S\setminus \tilde A$ is removable for $f$, see the Step 3 in the
proof of Theorem \ref{propa}. Therefore $f_k$ strongly converge
to $f$ outside of a proper analytic set $A\cup I_f$. We need now to
prove that $f$ is extendable to $\Delta^n$, \ie that $S$ is empty.
By Lemma \ref{area-bound2} below the areas of disks
$f_k(\Delta_{z^{'}})$ are bounded uniformly on $k$ and on $z^{'}\in
\Delta^{n-1}(1-\eps)$ for any fixed $\eps>0$, here
$\Delta_{z^{'}}\deff \{z^{'}\}\times \Delta$. Therefore the areas of
$f(\Delta_{z^{'}})$ are bounded to. Theorem 1.5 together with
Proposition 1.9 from \cite{Iv3} imply now that $f$ meromorphically
extends onto $\Delta^{n-1}(1-\eps)\times \Delta$. Therefore it
extends to $\Delta^n$. The condition (\sli of Definition
\ref{gamma-def} is fulfilled.

\smallskip Let $H$ be a
divisor in $X$. Then for every compact $K\comp H^n_{\eps}$ the
volumes of $f_k^*H\cap K$ counted with multiplicities are bounded.
By Oka-Riemenschneider theorem, see \cite{Rm}, the volumes of the
extensions of these divisors are bounded on compacts of $\Delta^n$
to. This verifies the condition (\slii  of Definition \ref{gamma-def}.
Proposition is proved.

\smallskip\qed

\newsect[PN]{Convergence of mappings with values in projective space}

Now let us examine  our notions of convergence on the example when the
target manifold is a complex projective space.

\newprg[PN.pres]{Meromorphic mappings to complex projective space}
Let a meromorphic mapping $f: U \to \pp^N$ be given. Without loss of generality
we suppose that the image of $f$ is not contained in a hyperplane. Then the
(complete) inverse image $f^{-1}(H)$ under $f$ of a hyperplane $H$
is a divisor in $U$. By $f^{-1}(\pp^n\setminus H)$ we shall understand
$U\setminus f^{-1}(H)$. Denote by $[w_0:w_1:...:w_N]$
the homogeneous coordinates  of  $\pp^N$. Let  $U_j=\{ w \in \pp^N:
w_j\neq 0 \}$ and let $\frac{w_0}{w_j},...,\frac{w_N}{w_j}$ be  affine
coordinates in  $U_j$. Set $D_j\deff f^{-1}(U_j)$, \ie $D_j = U\setminus
f^{-1}(H_j)$, where $H_j\deff \{w_j=0\}$. Since $U_0$ is isomorphic to
$\cc^N$ the restriction $f\arrowvert_{D_0}: D_0\longrightarrow U_0$
is given by holomorphic functions $\frac{w_1}{w_0}=f_1(z),...,\frac{w_N}{w_0}=f_N(z)$.
The coordinate change in  $\pp^N$ shows that
$f\arrowvert_{D_0 \cap D_j}: D_0 \cap D_j \longrightarrow \pp^N$
is given by  functions $\frac{w_1}{w_0}=\frac {1}
{f_j(z)},...,\frac{w_N}{w_0}=\frac {f_N(z)}{f_j(z)}$ which are holomorphic in
$D_j$. Therefore functions  $f_1,...,f_N$ are meromorphic on $D_0 \cup D_j$. This
proves that $f_1,...,f_N$ are meromorphic on $\bigcup_{j=0}^{N} D_j \subset U$.
We have that $U\setminus \bigcup_{j=0}^{N} D_j = \bigcap_{j=0}^Nf^{-1}(H_j)$,
\ie for every point from this set the image of every its neighborhood intersects
every $H_j$. Such point can be only an indeterminacy point of $f$. I.e.,
$U\setminus \bigcup_{j=0}^{N} D_j\subset I_f$. $I_f$ is analytic of codimension
$\ge 2$ and therefore by the theorem of Levi, see \cite{Lv} or \cite{Ha},
every $f_j$ meromorphically extends to $U$.

\smallskip If
$f_1\equiv...\equiv f_n\equiv 0$ then $f(U)\equiv 0 \in U_0$. If not, let
$f_1 \not \equiv 0$. One finds holomorphic functions $h_j$ et $g_j$  $ 0 \leq j \leq
N$ in a polydisk neighborhood $V$ of a given point $x\in U$, $g_j\neq 0$ such that
\[
 f_1= \frac{h_1}{g_1},...,f_N= \frac{h_N}{g_N}
\]
and therefore gets
\[
f:=\left[1: \frac{h_1}{g_1}:...:\frac{h_N}{g_N}\right] =
\left[\prod_{j=1}^{N}g_j : h_1  \prod_{j=2}^{N}g_j:...:h_N
\prod_{j=1}^{N-1}g_j\right].
\]

This proves that $f$ can be locally written in the form

\begin{equation}
 f(z):=[f_0(z):f_1(z):...:f_N(z)]
\end{equation}
as claimed.

\newprg[PN.weak]{Weak convergence of mappings with values in projective space}

Let us prove now the part (\slii of Theorem \ref{conv-pn} from the Introduction.
I.e.,

\begin{prop}
\label{weak-pn} A sequence of meromorphic mappings $f_k$ from a
complex manifold $U$ to $\pp^N$ converges weakly on compacts of $U$
if an only if for every point $z_0\in U$ there exists a neighborhood
$V\ni z_0$ and reduced representations $f_k=[f_k^0:...:f_k^N]$,
$f=[f^0:...:f^N]$ in $V$ such that for every $0\le j\le N$ $f_k^j$
converge to $f^j$ uniformly on $V$.
\end{prop}

\medskip\noindent $\Rightarrow $ Let $f_k\rightharpoonup f$, \ie $f_k$ converge
to $f$ weakly. Shrinking $U$ we suppose that  all $f_k$ and $f$ admit reduced
representations
\begin{equation}
\eqqno(reduced-fn)
f_k = [f_k^0:...:f_k^N]
\end{equation}
and
\begin{equation}
\eqqno(reduced-f)
f = [f^0:...:f^N]
\end{equation}
correspondingly. Up to making a linear coordinate change in $\pp^N$
we can suppose that $f[U]$ is not contained in any of coordinate
hyperplanes, \ie that $f^j\not\equiv 0$ for all $0\le j\le N$. Set

\[
Z^j = \{z\in U: f^j(z)=0\},
\]
and note that $\bigcap_{j=1}^NZ^j = I_f$. Since $f_k$ converge
on compacts in $U_j\deff U\setminus Z^j$ to $f$, see Remark
\ref{weak-rem}, we see, taking $j=0$, that

\begin{equation}
\eqqno(quot-conv1) \frac{f_k^j}{f_k^0} \rightrightarrows
\frac{f^j}{f^0}
\end{equation}
for all $j$ on compacts in $U_0$. Denote by  $Z_k^0$ the zero
divisors of $f_k^0$ and note that they leave every compact of $U_0$
as $n\to\infty$.

\begin{lem}
\label{div-conv} Divisors $Z_k^0$ converge to $Z^0$ in cycle space
topology.
\end{lem}

Let us prove this Lemma first. Fix a point $a\in Z^0\setminus Z^j$
(if  $Z^0\setminus Z^j$ is not empty) and take a relatively compact
neighborhood $V\ni a$ such that $\bar V\cap Z^j = \emptyset$. We
have that $f_k^0/f_k^j\rightrightarrows f^0/f^j$ on $\bar V$. The
Rouch\'e's theorem easily implies now that $Z_k^0\cap V$ converge to
$Z^0\cap V$ as currents.

\begin{rema} \rm
In fact the cycle space topology on the space of divisors coincides
with the topology of uniform convergence of defining them
holomorphic functions, see \cite{Stl}. And this immediately gives
the previous assertion.
\end{rema}

We conclude from here that $Z_k^0$ converge to $Z^0$ as cycles on
compacts in $U\setminus I_f$. But then by \cite{Ni}, Theorem II, we
obtain that they converge on the whole of $U$. Lemma \ref{div-conv}
is proved.

\medskip We continue the proof of the Theorem. Shrinking $U$ if
necessary we can suppose that $U$ is biholomorphic to $\Delta^n =
\Delta^{n-1} \times \Delta$ and $Z_k^0\cap U$ regularly covers
$\Delta^{n-1}$ for $k\gg 1$. Now each $Z_k^0$ can be written as the zero
set of a uniquely defined unitary polynomial $P_k$ from
$\calo_{\Delta^{n-1}}[z_n]$ and these $P_k$ uniformly converge to
$P$ - the defining polynomial for $Z^0$. After multiplying each
$[f_k^0:...:f_k^N]$ by the unit $P_k/f_k^0$ we get the reduced
representations

\[
f_k = [P_k:g_k^1...:g_k^N].
\]
The same with
\[
f = [P:g^1...:g^N].
\]
But now $P_k\rightrightarrows P$ and therefore from
\eqqref(quot-conv1), which reads now as
\begin{equation}
\eqqno(quot-conv2) \frac{g_k^j}{P_k} \rightrightarrows \frac{g^j}{P}
\end{equation}
on compacts in $U_0$, we get that for every $1\le j\le N$
$g_k^j\rightrightarrows g^j$ on compacts in $U_0 = U\setminus Z^0$.
But from the maximum principle it follows that
$g_k^j\rightrightarrows g^j$ on compacts in $U$.

\medskip\noindent $\Leftarrow $ For proving the inverse statement we
start with converging reduced representations \eqqref(reduced-fn) to
\eqqref(reduced-f), \ie $f_k^j \rightrightarrows f^j$ on $U$. Then
for every $0\le j\le N$  on every $U_j=U\setminus Z^j$ we get a
convergence on compacts
\[
\left(\frac{f_k^0}{f_k^j}, ..., \frac{f_k^N}{f_k^j}\right)
\rightrightarrows \left(\frac{f^0}{f^j}, ...,
\frac{f^N}{f^j}\right).
\]
And since the codimension of $I_f=\bigcap Z^j$ is at least two we
deduce the weak convergence of $f_k$ to $f$.

\smallskip\qed

\newprg[PN.strong]{Strong convergence and convergence of meromorphic functions}
Strong convergence of meromorphic maps into $\pp^N$ can be described in the
following way. First, if $f_k\to f$ then $f_k\rightharpoonup f$. Therefore
$[f_k^0:...:f_k^N]\rightrightarrows [f^0:...:f^N]$ for an appropriate reduced
representations. According to \eqqref(gr-vol2) the volume of the graph of $f_k$ is
\begin{equation}
\int\limits_{U\setminus I_{f_k}}\big(\omega_1 + f_k^{\ast }\omega_{FS}\big)^n =
\int\limits_{U\setminus I_f}\sum_{j=0}^nC_n^j\omega_1^j\wedge f_k^{\ast }\omega_{FS}^{n-j}.
\end{equation}
Since $f_k^*\omega_{FS} = dd^c\ln \norm{f_k}^2$ this is nothing but the non-pluripolar
Monge-Amp\`ere mass of $\ln \norm{f_k}^2$ as appeared in \eqqref(monge1). By
Proposition \ref{str-cycl} volumes of $\Gamma_{f_k}$ converge to the volume of
$\Gamma_f$, \ie

\begin{equation}
\eqqno(gr-vol3)
\int\limits_{U\setminus I_{f_k}}\omega_1^j\wedge \left(dd^c\ln \norm{f_k}^2\right)^{n-j}\to
\int\limits_{U\setminus I_f}\omega_1^j\wedge \left(dd^c\ln \norm{f}^2\right)^{n-j}
\end{equation}
for $0\le j<n$. In the case $U=\Delta^n$ this gives \eqqref(monge1). Vice versa, if
one has convergence of volumes the appearance of an exceptional component
is impossible and we conclude the part (\sliii of Theorem \ref{conv-pn}:

\begin{prop}
\label{strong-pn}
$f_k$ converge to $f$ strongly if and only if

\smallskip\sli the appropriate reduced representations converge uniformly;

\smallskip\slii for every $0\le j\le n-1$ one has \eqqref(gr-vol3).
\end{prop}

Now let us descend to the convergence of meromorphic functions.
Meromorphic functions on a complex manifold $U$ are exactly the meromorphic
mappings from $U$ to $\pp^1$. I.e., all our previous results and notions
are applicable to this case.

\begin{prop}
If a sequence $\{f_k\}$ of meromorphic functions converge weakly then it
converge strongly.
\end{prop}
\proof  Let $f$ be the weak limit of $f_k$.  We shall see in a
moment, see Corollary \ref{gma-bnd1} that volumes of graphs in this
case are uniformly bounded over compacts in $U$. Therefore after
going to a subsequence we get that the Hausdorff limit
$\hat\Gamma\deff\lim \Gamma_{f_k}$ is a purely $n$-dimensional
analytic subset of $U\times \pp^1$. We claim that $\lim\Gamma_{f_k}
= \Gamma_f$ in fact, \ie that there are no exceptional components.
If not take any irreducible component $\Gamma$ of this limit
different from $\Gamma_f$. Denote by $\gamma$ its projection to $U$.
$\gamma$ is a proper analytic set of codimension at least two $U$.
But then $\Gamma$ should be contained in $\gamma \times\pp^1$ and
the last analytic set is of dimension $\dim U - 1$. This is
impossible, because all components of $\lim \Gamma_{f_k}$ are of
pure dimension $\dim U$. Therefore $\gamma = \emptyset$ and
$\lim\Gamma_{f_k} = \Gamma_f$.

\smallskip\qed

\newprg[PN.gamma]{Gamma convergence in projective case}

In \cite{Fu} and subsequent papers of Fujimoto the following
type of convergence of meromorphic mappings with values in $\pp^N$
was considered, it was called the $m$-convergence (or {\slsf
meromorphic} convergence): $f_k$ $m$-converge to $f$ if there exist
reduced (admissible in the terminology of \cite{Fu}) representations
$f_k = [f_k^0:...:f_k^N]$ which converge uniformly on compacts to
$f=[f^0:...:f^N]$,  but the last {\slsf is not supposed} to be
reduced (\ie admissible), only not all $f^j$ are identically zero.
Let us prove the item (\sli of Theorem \ref{conv-pn}.

\begin{prop}
\label{gamma-pn}
When the target manifold $X$ is the complex projective space $\pp^N$ the $\Gamma$-convergence
of meromorphic mappings is equivalent to $m$-convergence in the sense of Fujimoto.
\end{prop}
\proof $\Rightarrow $ Suppose that $f_k\xrightarrow{\Gamma}f$. Let $\gamma$ be
the an analytic subset of $U$ such that our sequence converge strongly on compacts
of $U\setminus \gamma$. We add to $\gamma$ also the indeterminacies of the limit
$f$ and therefore $f_k$ will converge to $f$ on  $U\setminus \gamma$ in
compact open topology.  Let $f=[f^0:...:f^N]$ be some reduced representation of the limit map.
Making linear change of coordinates we can suppose, without loss of generality that $f^0\not\equiv 0$,
\ie that $f(U)\not\subset H_0$, where $H_0 = \{ Z_0 = 0\}$  in homogeneous coordinates $[Z_0:...:Z_N]$
of $\pp^N$. We have that $f_k^*H_0$ converge on compacts in $U$ in the cycle space topology
(after taking a subsequence).

\smallskip Take some $a\in \gamma$ and choose a chart $(V,j)$ adapted both to
$\gamma$ and $f^*H_0$ with coordinates $z_1,...,z_n$ around $a$ in such a way that $a=0$ and
$(\gamma\cup f^*H_0)\cap \left(\Delta^{n-1} \times\Delta\right)$ projects to $\Delta^{n-1}$
properly. Then $f_k^*H_0\cap V$ also projects to $\Delta^{n-1}$ properly for $k\gg 1$. After going to
a subsequence once more we can fix the degree $d$ of ramified coverings $f_k^*H_0\cap V\to
\Delta^{n-1}$ and write the corresponding polynomials $P_k\in \calo_{\Delta^{n-1}}[z_n]$ defining
$f_k^*H_0\cap V$. $P_k$ converge to some $P$ on (compacts of) $\Delta^{n-1}$. Let
$f_k=[f_k^0:...:f_k^N]$ be some reduced representations of $f_k$ on $V$. Notice that
$f_k^*H_0\cap V = \{f_k^0=0\}$. Divide each such representation by the unit $f_k^0/P_k$ and
get representations $f_k = [P_k:g_k^1:...:g_k^N]$ with converging first terms $P_k$. At the same
time $(g_k^1/P_k,..., g_k^N/P_k)$ represents $f_k$ in nonhomogeneous coordinates of the chart
$Z_0\not=0$ on $\pp^N$. Therefore $g_k^j/P_k$ converge to some $f^j$ on compacts of $V\setminus
(\gamma\cup f^*H_0)$. Therefore $g_k^j$ converge to $g^j\deff f^jP$ on compacts of
$V\setminus (\gamma\cup f^*H_0)$ to. By maximum principle they converge everywhere on $V$ to
the extension of $g^j$. We get that reduced representations $f_k = [P_k:g_k^1:...:g_k^N]$
converge term by term to a (may be non reduced) representation $[P:g^1:...:g^N]$ and this
can be only a representation of $f$.

\medskip\noindent $\Leftarrow $ Suppose now that $f_k$ $m$-converge to $f$. Again
change coordinates in $\pp^N$, if necessary, in such a way that $f(U)\not\subset H_0$.
Let $V$ be a neighborhood of some point $a\in U$. If $a\in f^*H_0$ then take
$(V,j)$ to be an adapted chart to this divisor. In any case take $V$ to be biholomorphic
to $\Delta^{n-1}\times \Delta$. Let $F_k = (f_k^0,...,f_k^N)$ be the lifts of $f_k$ to
$\cc^{N+1}$ in $V$ such that $F_k$ converge to the lift $F=(f^0,...,f^N)$ of
$f$. From here one gets immediately that $f_k^j/f_k^0$ converge to $f^j/f^0$ uniformly
on compacts of $V\setminus \{f^0=0\}$, \ie that our maps converge strongly outside
of a divisor.

\smallskip Now let $H = \{ P(Z_0,...,Z_N)=0\}$ be a divisor
such that $f(U)\not\subset H$. Using convergence of lifts $F_k=(f_k^0,...,f_k^N)$ to
$F=(f^0,...,f^N)$ one gets that $f_k(U)\not\subset H$ for $k\gg 1$. One has also that
$P(f_k^0,...,f_k^N)$ uniformly converge to $P(f^0,...,f^N)$ and this is equivalent
to the convergence of divisors.

\smallskip\qed

\begin{rema}\rm
The relation between weak/gamma convergence and $m$-convergence for
the case of $X=\pp^N$ was indicated without proof  in \cite{Iv2}.
\end{rema}

\newsect[NORM]{Bloch-Montel type normality criterion}

The aim of this section is to test the notion of weak convergence on
the Bloch-Montel type normality statement, \ie we are going to prove
here the Theorem \ref{montel} from the Introduction.

\newprg[NORM.prel]{Preliminaries} Before proceeding with the proof let us recall few basic
facts. We start with an extended version of Zalcman's lemma, see \cite{Me}:

\begin{lem}
\label{zalcman}
A family $\calf$ of holomorphic mappings from $\Delta^n$ to a compact Hermitian
manifold $(X,h)$ is not normal at $z_0\in \Delta^n$ if and only if there exist
sequences $z_k\to z_0$, $r_k\searrow 0$, $f_k\in \calf$ such that $f_k(z_k + r_kw)$
converge uniformly on compacts in $\cc^n$ to a non-constant entire mapping $f:\cc^n
\to X$ such that $\norm{df(w)}_h\le 2$ for all $w\in \cc^n$.
\end{lem}

This $f$ may well have rank one. We shall also need the following result from \cite{IS1},
which is a precise version of Gromov compactness theorem (we shall need it in the integrable
case only):

\begin{prop}
\label{stable-conv} Let $u_k:\Delta \to X$ be a sequence of
holomorphic maps into a disk-convex Hermitian manifold $(X,h)$ with
uniformly bounded areas, which uniformly converges on some annulus
$A_{1-\eps,1}$ adjacent to the boundary $\d\Delta$. Then $u_k$
converge to stable complex curve over $X$ after a reparametrization.
Moreover, the compact components of the limit are rational curves.
\end{prop}

For the notions of stable curve over $X$, convergence after a
reparametrization, as well as for the proof we refer to \cite{IS1}.
The obvious conclusion from this type of convergence is the
following:

\begin{corol}
\label{intersects}
If $u_k$ converge in stable sense to $u$ and $u(\Delta)$ intersects a divisor $H$ in $X$,
but us not contained in $H$, then all $u_k(\Delta)$ intersect $H$ for $k\gg 1$.
\end{corol}
\proof It was proved in \cite{IS2} (more details are given in \cite{IS1}) that for any $k\gg 1$
one can join $u_k$ with $u$ by a holomorphic one parameter family of stable maps, see Proposition
2.1.3 in \cite{IS2} for the exact statement. For us it is sufficient to understand that
there exists a normal complex surface $Y\xrightarrow{\pi} \Delta$ foliated over the disk $\Delta$
such that all fibers $Y_s\deff \pi^{-1}(s)$ are disks and a holomorphic mapping $\calu :
Y\to X$ such that $\calu|_{Y_0} = u$ and $\calu|_{Y_{s_0}} = u_k$ for some $s_0\in\Delta$ and
some $k$.

\begin{rema} \rm
The fact that this family can be contracted to a surface with normal points is proved in Lemma 2.2.6
in \cite{IS2}.
\end{rema}

Let $h$ be a defining holomorphic function of the divisor $H$ near the point of intersection
$u(\Delta)\cap H$. Then $h\circ \calu$ is holomorphic on $Y$ (for this one might need to take
disks of smaller radii) and is equal to zero at $0\in Y_0\subset Y$. At the same time it
cannot vanish on $\bigcup_{s\in \Delta}\d Y_s$ because $\calu|_{Y_s} (\d Y_s)$ is close to $u(\d\Delta)$
for all $s\in \Delta$.  Therefore the zero set of $h\circ \calu$ must intersect every $Y_s$.
And that means that $u_k(\Delta)$ intersects $H$.

\smallskip\qed

\smallskip Let us make one more remark.  Let $\omega_{FS}$ be
the Fubini-Study form on $\pp^N$. For a holomorphic map
$f:\bar\Delta \to \cc^N$ (we always suppose $f$ to be defined in a
neighborhood of the closure $\bar\Delta$), the area of $f(\Delta)$
with respect to the Fubini-Study form is

\begin{equation}
\eqqno(area-FS)
\area_{FS}f(\Delta) = \int\limits_{\Delta}f^*\omega_{FS}.
\end{equation}
Denote by $Z=(Z_0,...,Z_N)$ coordinates in $\cc^{N+1}$ and let $\pi:\cc^{N+1}\setminus \{0\}
\to \pp^N$ be the standard projection. Consider the following singular $(1,1)$-form on $\cc^{N+1}$
\begin{equation}
\eqqno(omega-0)
\omega_0 = dd^c\ln \norm{Z}^2.
\end{equation}

The following statement is a simple case of King's residue formula, but we shall give a simple
proof for the sake of completeness.

\begin{lem}
For a holomorphic lift $F=(f^0,...,f^N):\bar\Delta\to \cc^{N+1}$ of $f :\bar\Delta\to \pp^N$
(\ie $f=\pi\circ F$) such that $F|_{\d\Delta}$ doesn't vanishes one has

\begin{equation}
\eqqno(area-lift)
\area_{FS}f(\Delta) = \int\limits_{\d\Delta}d^c\ln \norm{F}^2 - N_F.
\end{equation}
Here $N_F$ is the number of zeroes of $F$ counted with
multiplicities.
\end{lem}
\proof  By the very definition of the Fubini-Study form one has
$\pi^*\omega_{FS} = \omega_0$. And therefore it is immediate to
check that in a neighborhood of a point $a\in \Delta$ such that
$F(a) \not=0$ one has that $f^*\omega_{FS} = F^*\omega_0$. As the
result

\begin{equation}
\eqqno(area-out-0)
\area_{FS}f(\Delta) = \int\limits_{\Delta}f^*\omega_{FS}=
\int\limits_{\Delta\setminus Z_F}F^*\omega_0,
\end{equation}
where $Z_F \deff \{z_1,...,z_k\}$ is the set of zeroes of $F$, \ie such $z_l$
that $f^j(z_l) = 0$ for all $j=0,...,N$. Let $n_l$ be the multiplicity of zero $z_l$.
Then $F(z) = (z-z_i)^{n_l}(g^0(z),...,g^N(z))$,
where at least one of $g^j$-s is not zero at $z_l$. We have that
\[
dd^c\ln \norm{F}^2 = n_l\delta_{z_l} + dd^c\ln\norm{G}^2,
\]
where $G(z) = (g^0(z),...,g^N(z))$. Therefore $dd^c\ln\norm{G}^2$ is an extension of $F^*\omega_0$ to
$z_l$. The rest obviously follows from the Stokes formula.

\smallskip\qed

Let us observe the following immediate corollary from this lemma.

\begin{corol}
\label{gma-bnd1}
Let $f_k : U\to \pp^N$ be a $\Gamma$-converging sequence of meromorphic mappings and
let $L$ be a divisor in $U$ such that $f_k$ converge uniformly on compacts of $U\setminus L$.
Let $V\cong \Delta^{n-1}\times \Delta$ be a scale adapted to $L$ and to the limit $M$
of $f_k^*H_0$, where $H_0 = [Z_0=0]$. Then the areas of the
analytic disks $f_k(\Delta_{z^{'}})$ are uniformly bounded in $z^{'}\in \Delta^{n-1}$ and
$k\in \nn$.
\end{corol}
\proof Let $(z^{'},z_n)$ be coordinates in $\Delta^{n-1}\times \Delta$. Denote by
$F_k = (f_k^0,...,f_k^N)$ lifts of $f_k$ to $\cc^{N+1}$. Consider
restrictions $f_k|_{\Delta_{z^{'}}}$. Due to the fact that our chart is adapted to
$M=\lim f_k^*H_0$ we have that $f_k^0$ doesn't vanishes on $\d \Delta_{z^{'}}$ for
$k\gg 1$ and, since it is also adapted to $L$ the lifts $F_k = (f_k^0,...,f_k^N)$ converge
in a neighborhood of $\d \Delta_{z^{'}}$. By \eqqref(area-lift) we have

\begin{equation}
\eqqno(area-Fk)
\area_{FS}f_k(\Delta_{z^{'}}) \le \int\limits_{\d \Delta_{z^{'}}}d^c\ln
\norm{F_k}^2\le c,
\end{equation}
\ie the areas are uniformly bounded for $z^{'}\in \Delta^{n-1}$ and all $k$.

\smallskip\qed

\begin{rema} \rm
\label{separ} For a family $\calf$ of meromorphic mappings from a
manifold $U$ to a projective manifold $X$ to be normal an obvious
necessary condition is that for any fixed hypersurface $H\subset X$
and any fixed compact $K\comp U$ the volumes counting with
multiplicities of intersections $f^*H\cap K$ should be uniformly
bounded for $f\in\calf$. It was proved by Fujimoto in \cite{Fu} that
this condition (in the case $X=\pp^N$ and $H_i$ are hyperplanes)
turns out to be also sufficient, but only for the {\slsf
meromorphic} (\ie $\Gamma$) normality.  We in this paper are interested in the
normality in the weak convergence sense (which is, that's to say,
stronger than meromorphic one). In that case there is one more necessary
condition. Take two hypersurfaces $H_0$ and $H_1$ in $X$. Let
$\{f^*H_0:f\in \calf\}$ and $\{f^*H_1:f\in \calf\}$ be the families
of their preimages by elements of our family $f\in \calf$. By
boundedness of volumes condition for every sequence $f_k^*H_i$,
$i=0,1$, some subsequences $f_{k_j}^*H_i$ converge to divisors $L_0$
and $L_1$. If there exist {\slsf coinciding} (without taking to account
the multiplicities) components $L_0^{'}$ and $L_1^{'}$ of  $L_0$ and
$L_1$ respectively, then $f_{k_j}$ cannot weakly converge in a
neighborhood of $L_0^{'} = L_1^{'}$. Indeed, since the limit $f$ is
a holomorphic map outside of $I_f$,  the preimages $f^*H_0$ and
$f^*H_1$ cannot have common components. But $L_0^{'}$ and $L_1^{'}$
are such components. Contradiction. This will be formalized in the
following definition.
\end{rema}

Let $\calf$ be a $\Gamma$-normal family in $\calm (U,X)$. Fix a
divisor $H$ in $X$. Remark that for every relatively compact $D\comp
U$ the intersections $f^*H\cap \bar D$ from a pre-compact family of
sets when $f$ is running over $\calf$. Therefore one can find a
finite collection of scales $\{E_{\alpha}\}$ such that every
$f^*H\cap \bar D$ can be covered by some of corresponding
$V_{\alpha}$ and the members of this covering are adapted to $f^*H$.
This collection $\{E_{\alpha}\}$ of scales depends on $D\comp U$,
but doesn't depend on $f\in \calf$ and, moreover, doesn't depend on
$H$ taken in some compact family of divisors, in our case this
family is $\{H_0,...,H_d\}$, \ie is finite.

\begin{defi}
\label{unif-separ}
We say that a family $\calf$ of meromorphic mappings from a complex
manifold $U$ to a complex manifold $X$ uniformly separates hypersurfaces
$H_0$ and $H_1$ from $X$ if for any $f\in \calf$ and any adapted for both
$f^*H_0$ and $f^*H_1$ scale $E_{\alpha}=(V_{\alpha},j_{\alpha},\Delta^{n-1}, \Delta)$
as above, the Hausdorff
distance between $f^*H_0\cap V_{\alpha}$ and $f^*H_1\cap V_{\alpha}$ for $f\in\calf$ is
bounded from below by a strictly positive constant.
\end{defi}

Hausdorff distance is taken here in the Euclidean metric of $\cc^n$.
A constant in question may well depend on divisors $H_0$, $H_1$ and
adapted chart $V_{\alpha}$, but it is supposed not to depend on $f\in\calf$.

\newprg[NORM.proof]{Proof of the normality criterion of Theorem \ref{montel}}
We are going to prove now Theorem \ref{montel} from the Introduction.
Recall that a relatively compact open subset $Y$ of a complex manifold
$X$ is said to be {\slsf hyperbolically imbedded} to $X$ if for any two
sequences  $\{ x_n\}$  and $\{ y_n\} $ in $Y$ converging to distinct
points $x \in \bar{Y}$ and $y \in \bar{Y}$  one has

\[
\lim\sup _{n\to \infty }k_Y(x_n,y_n)>0,
\]
where $k_Y$ is the Kobayashi pseudodistance of  $Y$. $Y\comp X$ is
said to be {\slsf locally  hyperbolically  complete} (l.h.c) if for
every $y\in \bar Y$ there exists a neighborhood $V_y\ni y$ such that
$V_y\cap Y$ is hyperbolically complete. For example every $Y\comp X$
of the form $X\setminus \{ \text{ divisor } \}$ is obviously l.c.h.
It was proved in \cite{Ki} that if $Y$ is hyperbolically imbedded
into  $X$ and is l.h.c. then  $Y$ is complete hyperbolic.

\smallskip These notions are connected to complex lines in $\bar Y$ by Theorem of
Zaidenberg, see \cite{Za}. By a complex line in $Y$ (or in $X$) one
understands an image of a non-constant holomorphic map  $u:\cc \to
Y$ (or $X$). Sometimes one requires that
$\norm{d_zu(\frac{\partial}{x})}_h \le 1$ for all $z\in \cc $, where
$h$ is some Hermitian metric on $X$. Complex line $u:\cc\to \bar
Y\comp X$ is called {\slsf limiting for $Y$} if there exists a
sequence of holomorphic mappings  $u_n :\Delta (R)\to Y$ converging
on compacts in $\cc$ to $u:\cc\to \bar Y$. Theorem of Zaidenberg
says now that: for a relatively compact l.c.h. domain $Y$ in a
complex manifold $X$ to be complete hyperbolic and hyperbolically
imbedded in $X$ it is necessary and sufficient that $Y$ doesn't
contain complex lines and doesn't admits limiting complex lines.

\medskip Now we turn to the proof. Let $\{f_k\}$ be a sequence from $\calf$, where
$\calf$ satisfies the assumptions of Theorem \ref{montel} from the
Introduction. $\{ H_i\}_{i=0}^{d}$ our set of divisors.


\medskip\noindent{\slsf Step 1. Convergence outside of a divisor.}
By Bishop's compactness theorem for every $i$ some subsequence from
$f_k^*H_i$ converges to a (may empty) hypersurface in $U$. Denote
this limit hypersurface as $L_i$. Set

\[
L \deff \bigcup_{i=0}^{d} L_i.
\]
In order not to complicate notations we will not introduce
subindexes when extracting subsequences.

\smallskip If  $L$ is empty then for every compact $K\comp U$
all $f_k$ with $k$ big enough send $K$ to $X \setminus\bigcup_{i=0}^d H_i$, the last is Stein.
In particular they are holomorphic in a neighborhood of $K$ and we can
use Zalcman's Lemma \ref{zalcman} together with Zaidenberg's characterization
to extract a converging subsequence.

\smallskip Therefore from now on we suppose  $L$ is nonempty. Take a point
$z_0\in U\setminus L$ and take a
relatively compact  neighborhood $V\ni z_0$ biholomorphic to a ball
such that $\bar V\cap L = \emptyset $. Then for $k$ big enough
$f_k(\bar V)\subset X\setminus \bigcup_{i=0}^{d} H_i$. This implies
that they all are holomorphic on $V$ and  we again can find a converging
subsequence on $V$ as before. Therefore some subsequence of $\{f_k\}$ (still
denoted as $\{f_k\}$) converge on compacts of $U\setminus L$ in the
usual sense of holomorphic mappings. Denote by $f$ its limit. $f$ is
a holomorphic map from $U\setminus L$ to $X$.

\smallskip\noindent{\slsf Step 2. Convergence across the divisor.}
Take a point $z_0\in L_{0}$, if $L_0$ is empty we can re-numerate
$L_i$-s. fix an imbedding $i:X\to \pp^N$ and let  $H$ be the
intersection of $X$ (\ie of $i(X)$) with hyperplane $\{ Z_0=0\}$ in
the standard homogeneous coordinates $[Z_0:...:Z_N]$ of $\pp^N$.
After going to a subsequence we have that $f_k^*H$ converge, denote
by $M$ the limit. Let $(V,j)$ be an adapted chart for  $L\cup M$
(and therefore also for $L_0$) at $z_0$ with the scale $E = (V,j,
\Delta^{n-1}, \Delta)$. Let $P_k[z_n]\in \calo_{\Delta^{n-1}}[z_n]$
be the defining unitary polynomial for $f_k^*H\cap V$. $P_k$
converges to the defining polynomial $P$ of $M\cap V$.

\smallskip  Let $f_k = [f_k^0:...:f_k^N]$ be reduced representations of $f_k$ on $V$
(we write $f_k$ for $i\circ f_k)$.
Then multiplying this representation by the unit $P_k/f_k^0$ we obtain a reduced
representation $f_k = [P_k:g_k^1:...:g_k^N]$. We have that $g_k^j/P_k$ converge on
compacts of $V\setminus (L\cup M)$. Therefore $g_k^j$ converge there to, denote by
$g^j$ its limit. We see that lifts $F_k = (P_k,g_k^1,...,g_k^N)$ converge to
$F \deff (P,g^1,...,g^N)$ on compacts of $V\setminus (L\cup M)$.  By maximum principle
they converge on $V$.

\smallskip In particular $f$ extends to a meromorphic mapping from $U$ to $X$.

\begin{rema} \rm
It is worth of noticing that at this stage we proved the
$\Gamma$-normality of our family. For the case $X=\pp^N$ with $H_i$
hyperplanes this was proved in \cite{Fu}. One more point worth of
noticing is that the extendibility of $f$ also follows from usual
complex  hyperbolic geometry, see \cite{Ko}.
\end{rema}

\smallskip\noindent{\slsf Step 3. Convergence outside of codimension two.}
 Changing indices of $H_i$, if necessary, we can suppose that our
family uniformly separates $H_0$ and $H_1$. Take a point $z_0\in
L_0\setminus \bigcup_{i\not=0}L_i$ such that $L_0$ in addition is
smooth at $z_0$.  Take an adapted scale $E=(V,j,\Delta^{n-1},
\Delta)$ for $L_0$ near $z_0$ which intersects $L$ only by the
smooth part of $L_0$ and, moreover, such that $j(L_0\cap V) = d\cdot
[\Delta^{n-1}\times \{0\}]$ for some multiplicity $d\ge 1$. Fix
coordinates $(z^{'},z_n)$ on $\Delta^{n-1}\times \Delta$. By
Corollary \ref{gma-bnd1} the areas of analytic disks
$f_k|_{\Delta_{z^{'}}}$ are uniformly bounded. Fix some $z^{'}\in
\Delta^{n-1}$ and take a subsequence $f_k$ such that
$f_k|_{\Delta_{z^{'}}}$ converge in stable topology to
$f|_{\Delta_{z^{'}}}$ plus a chain $C_{z^{'}}$ of rational curves.
By Corollary \ref{intersects} if $C_{z^{'}}$ intersects some $H_1$
with $i\not=0$ then $f_k|_{\Delta_{z^{'}}}(\Delta_{z^{'}})$
intersects $H_1$ to. But then $f_k^*H_1 \cap V$ is nonempty and
converge to $L_1\cap V$. This can be only $L_0\cap V$ with some
multiplicity, because $V$ was chosen in such a way that $L\cap V =
L_0\cap V$. The last violates the assumed uniform separability of
the pair $H_0, H_1$ by $\calf$. Therefore $C_{z^{'}}$ is empty. That
means that (some subsequence of) $f_k|_{\Delta_{z^{'}}}$ uniformly
on $\Delta_{z^{'}}$ converge to $f|_{\Delta_{z^{'}}}$. This implies
that the whole sequence $f_k$ restricted to $\Delta_{z^{'}}$
converge to $f$. Therefore $f_k$ converge to $f$ on $U\setminus
\Sing L$ in compact open topology as holomorphic mappings. This
proves the Theorem.

\smallskip\qed

\begin{rema} \rm
Theorem of Bloch, see also \cite{Gr}, states that $Y =
\pp^N\setminus \bigcup_{j=0}^{2N}H_i$ is hyperbolically imbedded to
$\pp^N$, where $H_i$ are hyperplanes in general position. Therefore
$Y = \pp^N\setminus \bigcup_{i=0}^{2N}H_i$ is an example for our
Theorem \ref{montel}.
\end{rema}

\newsect[VOL]{Behavior of volumes of graphs under weak and gamma convergence}

In this section we are concerned with the following question: let
meromorphic mappings $f_k:U\to X$ converge in some sense to a
meromorphic map $f$, what can be said about the behavior of volumes
of graphs of $f_k$ over compacts in $U$? If $f_k$ converge to $f$
strongly then, as it was proved in Theorem \ref{str-cycl}, for
every relative compact $V\comp U$ we have that
\begin{equation}
\eqqno(vol-conv)
\vol (\Gamma_{f_k|_V})\to \vol (\Gamma_{f|_V}).
\end{equation}
When $f_k$ converges only weakly one cannot, of course
expect anything like \eqqref(vol-conv). At most what one can expect
is that volumes of $\Gamma_{f_k}$ stay bounded over compacts in $U$
and converge to the volume of $\Gamma_f$ plus volumes of exceptional
components. I.e., the question is if for a weakly converging
sequence $\{f_k\}$ one has that for every relatively compact open
$V\comp U$ there exists a constant $C_V$ such that
\begin{equation}
\eqqno(vol-bnd) \vol (\Gamma_{f_k|_V}) \le C_V \qquad\text{ for all
} \quad k.
\end{equation}

This turns to be wrong in general, the following example was
communicated to us by A. Rashkovskii.

\newprg[VOL.monge]{Example of Rashkovskii}

\begin{exmp}\rm
\label{rash}
There exists a sequence $\eps_k\searrow 0$ such that {\slsf holomorphic}
mappings $f_k:\bb^3\to \pp^3$ defined as
\begin{equation}
\eqqno(rash1)
f_k : (z_1,z_2,z_3) \to [z_1:z_1 - \eps_k: z_2:z_3^k]
\end{equation}
converge weakly to $f(z) = [z_1:z_1:z_2:0]$ on compacts of the unit ball $\bb^3\subset \cc^3$,
but the volumes of graphs of $f_k$ over the ball $\bb^3(1/2)$ of radius $1/2$ diverge. In
fact
\begin{equation}
\eqqno(rash2)
\vol (\Gamma_{f_k})\cap (\bb^3(1/2)\times \pp^3) \ge k.
\end{equation}
\end{exmp}
Consider the following family of plurisubharmonic functions on the unit ball $\bb^3$ in $\cc^3$:
\begin{equation}
\eqqno(rash3)
u_{\eps , k}(z) = \ln (|z_1|^2 + |z_1 - \eps |^2 + |z_2|^2 + |z_3|^k), \quad \eps \in (0, 1/4).
\end{equation}
Note that every $u_{\eps ,k}$ is bounded in $\bb^3$ and its total MA
mass in $\bb^3(1/2)$ coincides with those of the function
\[
\tilde u_{\eps , k} \deff \max \{u_{\eps , k}, s_k\} \quad \text{where} \quad s_k = \min\{\ln
(|z_1|^2  + |z_2|^2 + |z_3|^k): z\in \sss^5(1/2)\}.
\]
Here $\sss^5(1/2) = \d\bb^3(1/2)$ is the sphere of radius $1/2$. This fact follows from the
Bedford-Taylor definition of the MA mass of a product of {\slsf bounded} psh functions,
see \cite{BT}: $dd^cu_1\wedge dd^cu_2 \deff dd^c(u_1dd^cu_2)$ and so on by induction. Here
the point is, of course, to prove that $dd^c(u_1dd^cu_2)$ is again a closed positive current.
Now one writes
\[
\ma_{\bb^3(1/2)}(u_{\eps , k}) = \int\limits_{\bb^3(1/2)}(dd^cu_{\eps , k})^3 =
\int\limits_{\bb^3(1/2)}dd^cu_{\eps , k}\wedge (dd^cu_{\eps , k})^2 =
\int\limits_{\d\bb^3(1/2)}d^cu_{\eps , k}\wedge (dd^cu_{\eps , k})^2 =
\]
\[
= \int\limits_{\d\bb^3(1/2)}d^c\tilde u_{\eps , k}\wedge (dd^c\tilde u_{\eps , k})^2 =
\int\limits_{\bb^3(1/2)}(dd^c\tilde u_{\eps , k})^3 =\ma_{\bb^3(1/2)}(\tilde u_{\eps , k})
\]
because $u_{\eps , k} = \tilde u_{\eps , k}$ on the sphere $\sss^5(1/2)$. Since
$\tilde u_{\eps , k}$ converge uniformly to $\tilde u_k = \max\{ \ln (2|z_1|^2 + |z_2|^2 + |z_3|^k),
s_k\}$ as $\eps \searrow 0$ and $\ma_{\bb^3(1/2)}(\tilde u_k) = \ma_{\bb^3(1/2)}(u_k) = 4k$,
where $u_k = \ln (2|z_1|^2 + |z_2|^2 + |z_3|^k)$, we shall have that for $\eps_k$
small enough $\ma_{\bb^3(1/2)}(u_{\eps_k,k})\ge k$. This finishes the proof.

\begin{rema} \rm
\label{kisel} Examples of psh functions with polar singularities and
unbounded non-polar MA mass where constructed first by Shiffman and
Taylor, see \cite{Si1}, and especially simple one by Kiselman, see
\cite{Ks}: $u(z_1,...,z_n) = (1-|z_n|^2)(-\ln \norm{z^{'}}^2)^{1/2}$
for $z^{'} = (z_1,...,z_{n-1})$. Taking any of these examples and
smoothing it by convolutions one gets a decreasing sequences of psh
functions converging outside of an analytic set (on any codimension)
to a  psh function, smooth outside of this set with unbounded
non-polar MA mass. The remarkable feature of the example of
Rashkovskii, just described, is that functions in this example have
a {\slsf geometric meaning}, their $dd^c$-s are pullbacks of
Fubini-Study form by a meromorphic mappings to the complex
projective space, \ie the sum of their non-polar MA masses are the
volumes of the corresponding graphs.
\end{rema}

\newprg[VOL.dim2]{Case of dimensions one and  two}

If $\{f_k\}$ is a $\Gamma$-converging sequence of meromorphic mappings
with values in one dimensional complex manifold then it is easy to
see that the volumes of graphs of $f_k$-s are locally bounded over
compacts in the source. Indeed, a one dimensional manifold  $X$
either properly imbeds to $\cc^n$ (when $X$ is noncompact) or is
projective and therefore imbeds to $\pp^n$. In both cases by Theorem
\ref{conv-pn} we have convergence of reduced representations to a,
may be nonreduced representation of the limit. Inequality
\eqqref(area-Fk) implies that in an appropriately chosen local
coordinates $(z^{'},z_n)$ one has
\[
\vol (\Gamma_{f_k|_{\Delta^n}}) = \int\limits_{\Delta^n}(dd^c||z||^2)^n +
\int\limits_{\Delta^n}(dd^c||z||^2)^{n-1}\wedge f_k^*\omega_{FS} \le
\int\limits_{\Delta^n}(dd^c||z||^2)^n +
\]
\[
+ \int\limits_{\Delta^{n-1}}(dd^c||z||^2)^{n-1} \int\limits_{\d \Delta_{z^{'}}}d^c\ln
\norm{F_k}^2\le \const.
\]

Next, if the dimension $n$ of the source $U$ is $2$ the boundedness of
volumes of graphs of a weakly converging sequence is automatic. This
can be seen at least in two ways.
First, in projective case this readily follows from the following formula of King,
see \cite{Kg}:
\begin{equation}
\eqqno(king)
d\left[d^c\ln (\norm{f}^2)\wedge dd^c\ln (\norm{f}^2)\right] = \chi_{U\setminus I_f}
\left[\left(dd^c\ln (\norm{f}^2)\right)^2\right] -  \sum_jn_j\left[Z_j\right],
\end{equation}
provided $I_f$ has pure codimension two. $Z_j$ are irreducible
components (branches) of the indeterminacy set $I_f$ of $f$. If
it has branches of higher codimension then around these branches a
higher order non-pluripolar masses can be expressed in a similar
way. Now if $f_k$ weakly converge to $f$ formula \eqqref(king)
immediately gives a uniform bound of corresponding MA masses (even
together with that concentrated on pluripolar sets $I_{f_k}$). If
$n=2$ then that's all we need.

\smallskip Second, using Skoda  potentials,  or Green functions, as
it was done in \cite{Iv2}
Theorem 2, one can bound non-pluripolar Monge-Amp\`ere masses of
order two also in the case of weakly converging sequence with values
in disk-convex K\"ahler $X$. This observation implies that if $X$ is
disk-convex K\"ahler and $\dim U =  2$ then the volumes of graphs of
weakly converging sequences of meromorphic mappings $U\to X$ are
uniformly bounded over compacts in $U$.

\smallskip
Moreover, it was proved in \cite{Ne} that volumes of weakly
converging sequence are bounded also in the case when $X$ is any
compact complex surface. The proof uses Ka\"ahler case separately
and then the fact that a non-K\"ahler surface has only finitely many
rational curves.

\begin{rema} \rm
Let us remark that there is one more important
case when the volumes of graphs of weakly (even $\Gamma$) converging
sequence necessarily stay bounded: namely when  $\{f_k\}$ is a  $
\Gamma$-converging sequence of meromorphic mappings between projective
manifolds $X$ and $Y$. Indeed the volumes of graphs $\Gamma_{f_k}$ are
uniformly bounded as it is straightforward from Besout theorem.
\end{rema}

\newsect[RAT-CON]{Rational connectivity of the exceptional components of the limit}

\newprg[RAT-CON.curves]{Chains of rational curves}

Recall that a rational curve $C$ in a complex manifold $X$ is an
image of $\pp^1$ in $X$ under a non-constant holomorphic map
$h:\pp^1\to X$. A chain of rational  curves is a {\slsf connected}
union $C=\bigcup_jC_j$ of finitely many rational curves.

\begin{defi}
A closed subset  $\Gamma \subset X$ we call rationally connected if for very two points
$p \not= q$ in $\Gamma $ there exists a chain of rational curves $C\subset \Gamma $
such that $p,q\in C$.
\end{defi}

One says also that $C$ {\slsf connects} $p$ with $q$. If $\Gamma $ is a complex manifold then
this property is equivalent to the either of the following two ones:

\begin{itemize}
\smallskip\item Every two points in $X$ can be connected by a single rational curve.

\smallskip\item For any finite set of points  $F\subset X$ there exists a rational curve
$C\supset F$.
\end{itemize}

\smallskip We refer to  \cite{Ar} for these facts. Now let us turn to the proof of Theorem
\ref{rat-con0} from the Introduction. It consists from the two
following lemmas. Let $f_k$ be a weakly or, gamma-converging
sequence of meromorphic mappings and $f$ denotes their limit. Let
$\hat\Gamma$ be the Hausdorff limit of the graphs, $\Gamma =
\overline{\hat\Gamma\setminus \Gamma_f}$ the corresponding bubble.
Set $\gamma \deff \pr_1(\Gamma)$. It is at most a divisor in
$\Gamma$-case and has codimension $\ge 2$ in the weak case. Let
$V\cong \Delta^{n-1}\times \Delta$ be a scale adapted to $\gamma $
in the sense that $(\bar\Delta^{n-1}\times\Delta)\cap \gamma
=\emptyset$.

\begin{lem}
\label{area-bound1}

Suppose that there exists a dense subset $S\subset \Delta^{n-1}$ such that the areas
of the analytic disks $\Gamma_{f_k|_{\Delta_{z^{'}}}}$ are uniformly bounded in $z^{'}\in S$
and $k\in \nn$ then for every point $a\in \gamma$ the fiber $\Gamma_a \deff \pr_2(\pr_1^{-1}(a))$
is rationally connected.
\end{lem}
\proof Here writing $f_k|_{\Delta_{z^{'}}}$ we mean the restriction
of $f_k$ to the disk $\Delta_{z'}\deff \{z'\}\times \Delta$. Fix a
point $a\in \gamma$ and some $a^1,a^2\in \Gamma_a$. Suppose
$a^1\not=a^2$, otherwise there is nothing to prove. We need to prove
that there exists a chain of rational curves in $\Gamma_a$
connecting  $a^1$ with $a^2$. Since $\Gamma_{f_k}$ converge to
$\hat\Gamma\supset \Gamma_a$ there exist $a_k^1\to a$ and $a_k^2\to
a$ such that $f_k(a_k^1)\to a^1$ and $f_k(a_k^2)\to a^2$. Perturbing
slightly we can take  such $a_k^i$ to be regular (\ie not
indeterminacy)  points of $f_k$ for $i=1,2$. Take a scale adapted to
$\gamma$ near $a$ in the sense that $\gamma\cap (\Delta^{n-1}\times
\d\Delta) = \emptyset$. Denote by $(z',z_n) = (z_1,...,z_{n-1},
z_n)$ the corresponding coordinates and assume without loss of
generality that $a=0$.

\smallskip
Let $b_k^1\to 0^{'}$ and  $b_k^2\to 0^{'}$ in $\Delta^{n-1}$  be such that $a_k^1\in
\Delta_{b_k^1}$ and $a_k^2\in \Delta_{b_k^2}$. Taking again $a_k^i$ sufficiently
general we can arrange that $b_k^i\in S$ and disks $\Delta_{b_k^1}$ and
$\Delta_{b_k^2}$ converge to the disk $\Delta_{0^{'}}$. After taking a
subsequence we get that graphs in question converge to the graph $\Gamma_{f|_{\Delta_{0^{'}}}}
\cup C^i$, where $C^i\subset \{a\}\times X$ are chains of rational curves. Both
these chains contain the point $f|_{\Delta_{0^{'}}}(a)$. Therefore $C\deff C^1\cup C^2$ is
connected. At the same time by construction $C^i\ni a^i$. Lemma is proved.

\smallskip\qed

\newprg[RAT-CON.proof]{Proof of Theorem \ref{rat-con0}}
Let us first consider the case of $\Gamma$-converging sequence of meromorphic mappings with
values in projective $X$.  Corollary \ref{gma-bnd1} gives us the required boundedness
of ares of analytic disks which makes possible to apply Lemma \ref{area-bound1} just 
proved. This gives us the statement of Theorem \ref{rat-con0} for $\Gamma$-converging 
sequences of meromorphic mappings with values in {\slsf projective} manifolds.

\smallskip To treat the case of Gauduchon target manifolds we shall need one more lemma.

\begin{lem}
\label{area-bound2} Let ${\cal F}$ be a family of meromorphic mappings
from $\Delta^n$ to a disk-convex manifold $X,$ which admits a
pluriclosed metric form. Suppose that for some  $0 < \epsilon < 1,$
the family ${\cal F}$ is holomorphic and equicontinuous on the
Hartogs figure $H^n_{\eps}$. Then for every $0 < r < 1$ the areas of
graphs $\Gamma_{f_{z'}}$ of restrictions  $f_{z'}\deff
f|_{\Delta_{z'}(r)}$ of $f\in\calf$ to the disks
$\Delta_{z'}(r)\deff \{z'\}\times \Delta_r$ are uniformly bounded in
$z'\in\Delta^{n-1}_r$ and $f\in\calf$.
\end{lem}
\proof For $f: \Delta^n \longrightarrow X $ a meromorphic map, we
denote by $I_ {f} \subset  \Delta^n$ the  set of points of
indeterminacy  of $f$. Since we suppose that all $f\in \calf$ are
holomorphic on $H^n_{\eps}$ the sets $I_f$ do not intersect
$\Delta^{n-1}\times A_{1-\eps ,1}$. Consider currents $T_{f} =
f^\ast\omega$ on $\Delta^n$, where $\omega$ is a pluriclosed metric
form on $X$. Write
\[
T_f =\frac{i}{2 } t^{\alpha\bar{\beta}}_f dz_\alpha \wedge d\bar{z_\beta},
\]
where $t^{\alpha\bar{\beta}}_f$ are distributions on $\Delta^n$ (in fact measures),
smooth on $\Delta^n\setminus I_{f}\supset H^n_{\eps}$. Fix $1-\eps < r < r_1 < 1$
and consider on $\Delta^{n-1}\setminus \pi (I_f)$ (where $\pi : \Delta^n\to\Delta^{n-1}$
is the canonical projection onto the first factor) the {\slsf area functions} $a_{f}$ given by
\begin{equation}
\eqqno(a-f)
a_{f}(z')= \area f_{z'}(\Delta_{r_1}) = \int_{\Delta_{z'}(r_1)}T_{f} =
\frac{i}{2}\int_{\Delta_{z'}(r_1)} t^{n\bar{n}}_{f} d{z_n} \wedge d\bar{z_n}.
\end{equation}
Functions $a_{f}$ are well-defined and smooth on $\Delta^{n-1}\setminus \pi (I_f)$.

The proof of Proposition will be done in two steps.

\smallskip\noindent{\slsf Step 1. }{\it Distributions $t^{n\bar n}_f $ are
locally integrable in $\Delta^n$.}
Note that forms  $T_{f}$ are smooth on $H^n_{\eps}$ and the family
$\{T_{f}: f \in{\cal F}\}$ is equicontinuous there.  The condition that $dd^cT_f=0$ implies,
in particular, that for all
$1\le k,l\le n-1$ one has
\begin{equation}
\eqqno(pl-clos)
\frac {\partial^2 t^{n\bar{n}}_{f}}{\partial{z_k}\partial {\bar{z_l}}}+ \frac{
\partial^2 t^{k\bar{l}}_{f}}{\partial{z_n}\partial{\bar{z_n}}}-\frac{ \partial^2
t^{k\bar{n}}_{f}}{\partial{z_n}\partial{\bar{z_l}}}-\frac{\partial^2
t^{n\bar{l}}_{f}}{\partial{z_k}\partial{\bar{z_n}}}=0.
\end{equation}
From \eqqref(pl-clos) we get that on $\Delta^{n-1}\setminus \pi (I_f)$:
\begin{equation}
\eqqno(ddc-a)
dd^c a_{f} =\left(\frac{i}{2}\right)^2 \sum_{k,l=1}^{n-1} \left(\int\limits_{\Delta_{z'}(r_1)}
\frac {\partial^2 t^{n\bar{n}}_{f}}{\partial{z_k}\partial {\bar{z_l}}} d{z_n} \wedge d\bar{z_n}\right)
dz_k\wedge d\bar z_l =
\end{equation}
\[
= \left(\frac{i}{2}\right)^2 \sum_{k,l = 1}^{n-1}\int\limits_{\Delta_{z'}(r_1)}
\left( \frac{ \partial^2 t^{k\bar{n}}_{f}}{\partial{z_n}\partial{\bar{z_l}}}+
\frac{ \partial^2 t^{n\bar{l}}_{f}}{\partial{z_k}\partial{\bar{z_n}}} -
\frac{ \partial^2 t^{k\bar{l}}_{f}}{\partial{z_n}\partial{\bar{z_n}}}\right)
d{z_n} \wedge d\bar{z_n}\cdot dz_k\wedge d\bar z_l =
\]
\[
= \left(\frac{i}{2}\right)^2\sum_{k,l = 1}^{n-1}\left(\int\limits_{\partial\Delta_{z'}(r_1)} \frac{ \partial
t^{k\bar{n}}_{f}}{\partial{\bar{z_k}}} d\bar{z_n} +  \int\limits_{\partial\Delta_{z'}(r_1)} \frac{ \partial
t^{n\bar{l}}_{f}}{\partial{{z_k}}} d{z_n} - \int\limits_{\partial\Delta_{z'}(r_1)}
\frac{ \partial t^{k\bar{l}}_{f}}{\partial{\bar{z_n}}} d\bar{z_n}\right) dz_k\wedge d\bar z_l =:\varphi_{f}.
\]

Forms $\varphi_{f}$ are smooth in the whole unit polydisk $\Delta ^{n-1}$ and equicontinuous there
because forms $T_{f}$ are smooth in  $\Delta^{n-1}\times A_{1-\eps, 1} \subset H^n_{\eps}$ and
equicontinuous there. Let us find a smooth and equicontinuous family on $\Delta^n_r$ of solutions
$\psi_{f}$ of
\begin{equation}
\eqqno(ddc-1)
dd^c \psi_{f}=\varphi_{f}.
\end{equation}
Set
\begin{equation}
\eqqno(h-f)
h_{f}:= a_{f}-\psi_{f}.
\end{equation}
Since $a_{f}$ is positive on $\Delta^{n-1}\setminus \pi (I_f)$ and
$\psi_{f}$ is smooth on $\Delta^{n-1}$ we see that  $h_{f}$ is bounded on $\Delta^{n-1}$ from below.
Also $dd^c h_{f}=0$ on $\Delta^{n-1}\setminus \pi (I_f)$ and therefore $h_{f}$ extends to a
plurisuperharmonic function on $\Delta^{n-1}$. This implies that $h_{f}\in
L^1_{loc}(\Delta^{n-1})$ see \cite{Ho}. It follows that $a_f$ and $t^{n\bar{n}}_{f}$ are locally
integrable. Step 1 is proved.

\smallskip\noindent{\slsf Step 2.}{\it Under the hypotheses of Lemma \ref{area-bound2}
functions $a_f$ defined by \eqqref(a-f) are smooth on $\Delta^{n-1}$ and for every fixed $r <1$
the family $\left\lbrace a_f\right\rbrace_{f\in \cal F}$ is  equicontinuous on
$\bar \Delta^{n-1}_r$.} Function $h_f$ given by  \eqqref(h-f) is plurisuperharmonic
in $\Delta^{n-1}$ and pluriharmonic on $\Delta^{n-1}\setminus \pi (I_f)$. Therefore by Siu's lower
semicontinuity of the level sets of Lelong numbers we have
\begin{equation}
\eqqno(ddc-h-f)
dd^ch_f = - \sum_{A \text{ irr.comp. of } \pi(I_f)} c_A(f)[A],
\end{equation}
where $c_A(f)\ge 0$ and $[A]$ denotes the current of integration over the irreducible component
$A$ of $\pi (I_f)$ of codimension one.

\begin{rema} \rm
Note that through components of higher codimension a pluriharmonic function $h_f$ extends (as
a pluriharmonic function). Therefore in \eqqref(ddc-h-f) the sum is taken over the components
of codimension one only.
\end{rema}
We need to prove that $c_A(f)=0$. From \eqqref(h-f) we get
\begin{equation}
\eqqno(ddc-af1)
dd^c a_f = dd^c \psi_f - \sum_{A \text{ irr.comp. of } \pi(I_f)} c_A(f)[A],
\end{equation}
where $dd^c$ from $a_f$ is taken in the sense of distributions (as from $L^1_{loc}$-function).
Let $\{h_A\}$ be equations of $A$. By Poincar\'e formula, see \cite{GK}, $[A]=dd^c\ln |h_A|^2$
and therefore \eqqref(ddc-af1) writes as
\begin{equation}
dd^c a_f = dd^c \psi_f - \sum_{A \text{ irr.comp. of } \pi(I_f)} c_A(f)dd^c\ln |h_A|^2
\end{equation}
Take an one dimensional disk $\Delta$ in $\Delta^{n-1}$ which intersects $\pi (I_f)$ transversely
at points $\{z_j\}$. Then \eqqref(ddc-af1) gives for restrictions of $a_f$ and $\phi_f$ to
$\Delta$ (and we shall denote them by the same letters) the following
\begin{equation}
\eqqno(ddc-af2)
\Delta a_f= \Delta \psi_f -\sum_{z_j \in \pi(A_f)} c_j(f)\delta_{z_j(f)}.
\end{equation}
Fix $\delta > 0$ such that  $\Delta(\delta,z_j)$  are pairwise disjoint. Let $\phi$ be a test
function on $\Delta$ with support in $\Delta (\delta,z_j)$ for some fixed $j$. The coordinate on
$\Delta$ denote as $z_1$.

\smallskip
Set
\[
a^{\epsilon}_f(z_1)=\frac{i}{2} \int\limits_{\Delta_{z_1}(r_1)}
t^{n\bar{n}}_{f,\epsilon} d{z_n} \wedge d\bar{z_n},
\]
where $t^{n\bar{n}}_{f,\epsilon}$ is the smoothing of $t^{n\bar{n}}_{f}$ by convolution.
Since $t^{n\bar{n}}_{f,\epsilon} \to t^{n\bar{n}}_{f}$ in  $L^1_{loc}$ we get by Fubini
Theorem  that $ a^{\epsilon}_f \to a_f$ in $L^1_{loc}$. Therefore using \eqqref(ddc-a)
for dimension two we obtain
\[
< \Delta a^{\epsilon}_f,\phi > = \frac{i}{2}\int\limits_{\Delta (\delta, z_j)} \phi (z_1) \left(
\int\limits_{\Delta_{z_1}(r_1)} \frac{\partial^2 t^{n\bar{n}}_{f,\epsilon}}{\partial
z_1 \partial \bar z_1} d{z_n} \wedge d\bar{z_n} \right) dz_1\wedge d\bar z_1 =
\]
\[
= \frac{i}{2 }\int\limits_{\Delta (\delta, z_j)} \phi
(z_1) \left( \frac{i}{2} \int\limits_{\partial\Delta_{z_1}(r_1)} \frac{ \partial
t^{1\bar{2}}_{f,\epsilon}}{\partial{\bar{z_1}}} d\bar{z_2} + \frac{i}{2}
\int\limits_{\partial\Delta_{z_1}(r_1)} \frac { \partial
t^{2\bar{1}}_{f,\epsilon}}{\partial{{z_1}}} d{z_2} \right)dz_1\wedge d\bar z_1 -
\]
\[
- \frac{i}{2 }\int\limits_{\Delta (\delta, z_j)} \phi (z_1) \left( \int\limits_{\partial
\Delta_{z_1}(r_1)} \frac{ \partial t^{1\bar{1}}_{f,\epsilon}}{\partial{\bar{z_2}}}
d\bar{z_2}\right) dz_1\wedge d\bar z_1 \longrightarrow < \varphi_{f},\phi >
\]
as $\epsilon \longrightarrow 0$. Therefore, $\Delta a_f= \varphi_f $
in $\Delta$ in the sense of distributions. By regularity of the
Laplacian $a_f \in \cal C^{\infty}$ on $\Delta$ and therefore
$c_A(f)=0$ for all $A$ and all $f$. Therefore $a_f$ are smooth on
$\Delta^{n-1}$ and $a_f = \psi_f + h_f$ there. $\psi_f$-s are
equicontinuous and $h_f$ are pluriharmonic everywhere and uniformly
bounded from below. Moreover $a_f$ are equicontinuous on
$\Delta^{n-1}_{\eps}$ by assumption. Therefore $h_f$ are
equicontinuous on $\Delta^{n-1}_{\eps}$. This implies equicontinuity
of $h_f$ on compacts of $\Delta^{n-1}$, and therefore the
equicontinuity of $a_f$. Step 2 and therefore our Lemma are proved.

\smallskip\qed

\smallskip Lemmas \ref{area-bound1} and \ref{area-bound2} obviously
imply the Theorem \ref{rat-con0} from the Introduction for the case
of weakly converging sequences of meromorphic mappings with values 
in disk-convex Gauduchon manifolds.

\newsect[FATOU]{Fatou components}

\newprg[FATOU.d-2]{Case of dimension two and Fatou sets}

First let us prove two lemmas.

\begin{lem}
\label{rat-con3} Suppose that a weakly converging sequence $\{f_k\}$
of meromorphic mappings from a two-dimensional domain $U$ to a
compact complex surface $X$ doesn't converge strongly. Then $X$ is
bimeromorphic to $\pp^2$.
\end{lem}
\proof Indeed, in that case there exists a point $a\in U$ and a
neighborhood $V\ni a$ such that $f_k$ converge uniformly on compacts
of  $V\setminus \{a\}$ but $\Gamma_{f_k}$ do not converge to
$\Gamma_f$, where $f:U\to X$ is the limit map. $\vol (\Gamma_{f_k})$
are uniformly bounded. Indeed, for K\"ahler $X$ it was proved in
\cite{Iv2} using Skoda's potentials. In \cite{Ne} its was proved for
non-K\"ahler $X$ using that fact that such $X$ can contain only
finitely many rational curves as well as existence of certain $dd^c$-exact
$(2,2)$-forms.

\smallskip Therefore we see that the limit $\hat\Gamma = \lim \Gamma_{f_k}$
contains  $\Gamma_f$ plus $\{a\}\times X$ (with some
multiplicity). But this is a bubble and therefore $X$ is rationally
connected by Theorem \ref{rat-con0}. From the classification of surfaces,
see \cite{BPV}, we know that such $X$ must be bimeromorphic to $\pp^2$.

\smallskip\qed

For a meromorphic map $f:U\to X$ denote by
\[
D_f\deff \pr_1\left(\{(z,x)\in U\times X: \dim_{(z,x)}\pr_2^{-1}(x)
\ge 1 \}\right)
\]
the set of degeneration of $f$. $f:U\to X$ is degenerate if $D_f=U$.

\begin{lem}
\label{degen} Let $f:X\to X$ be a non-degenerate meromorphic self-map
of a compact complex surface $X$ and let $z\in \Phi_s$ (resp.
$\Phi_w$). Then for every $l\ge 1$ one has
\begin{equation}
\eqqno(fat-iter)
f^l[z]\setminus f^l|_{D_{f^l}}(D_{f^l})\subset \Phi_s \quad\text{
(resp. } \Phi_w \text{)}.
\end{equation}
\end{lem}
\proof Take some $a\in f^l[z]\setminus f^l|_{D_{f^l}}(D_{f^l})$.
Since $f^l$ doesn't contract any curve to $a$ there exist
neighborhoods $V\ni z$ and $U\ni a$ such that $\pr_2 :(V\times
U)\cap\Gamma_{f^l} \to U$ is proper. That means that
$(f^l)^{-1}:U\to V$ is well defined as a multivalued holomorphic
map. Now let $\{f^{n_k}\}\subset \{f^n\}$ be a subsequence. By
assumption from the sequence $\{f^{n_k+l}\}$ we can subtract a
strongly/weakly converging on $V$ subsequence $\{f^{n_{k_j}+l}\}$.
That means that $\{f^{n_{k_j}}=f^{n_{k_j}+l}\circ f^{-l}\}$ will
converge in an appropriate sense on $U$.

\smallskip\qed

\smallskip Let us turn to the proof of Corollary \ref{fatou} from
the Introduction. Since every compact complex surface admits a
$dd^c$-closed metric form Theorem \ref{propa} applies in our case
and gives local pseudoconvexity of the weak Fatou set $\Phi_w$.
Suppose now that $\Phi_s \not= \Phi_w$.

\smallskip\noindent ${\sf a)}$ By Lemma \ref{rat-con3} $X \backsimeq \pp^2$.

\smallskip\noindent ${\sf b)}$ There exists a point $p\in X$, a ball
$B$ centered at $p$, a subsequence of iterates $\{ f^{n_k}\} $,
which uniformly converges on compacts of $\bar B \setminus \{ p\} $
to a meromorphic map $f_{\infty }:\bar B \to X$, holomorphic on
$B\setminus \{p\}$, but not converges strongly on any neighborhood of
$p$. In particular this means that $p\in I(f_{\infty })$ by Rouch\'e
Principle of \cite{Iv2} and, moreover, $C=f_{\infty }[p]$ is a chain
of rational curves $\bigcup_{i=1}^NC_i$. As it was said $\vol
(\Gamma_{f^{n_k}})$ are uniformly bounded on $\bar B $. So $\Gamma_{f^{n_k}}$
converge (after going to a subsequence) in cycle topology to
$\Gamma_{f_{\infty }}\cup d(\{ p\} \times X)$ for some integer $d\ge
1$. In particular $f$ cannot be degenerated in this case. 
Take a point $q\in X\setminus C$. Then for $k\gg 1$ we have that
$q\in f^{n_k}(B \setminus \{ p\} )$. If moreover $q\not\in
f^{n_k}(D(f^{n_k}))$ then $q\in \Phi_w$. But $\bigcup_k
f^{n_k}(D(f^{n_k}))$ is at most countable set of points and $\Phi_w$
is Levi-pseudoconvex. So $\Phi_w\supset X\setminus C$. Again from
pseudoconvexity of $\Phi_w$ it follows that if $\Phi_w$ intersects
some irreducible component of $C$ then it contains this component
minus the rest of $C$. I.e., $\Phi_w = \pp^2\setminus \{\text{ some
components of } C\}$.

\smallskip\noindent ${\sf c)}$ Take a point $(p,x)\in \{ p \}\times
X $ such that $x\in C$. Suppose that $\Gamma_{f_{\infty }}\cap (X
\times \{x\})$ has $(p,x)$ as isolated point. Then we can find
neighborhoods $W\ni p$ and $V\ni x$ such that $(\d W\times \bar V
)\cap \Gamma_{f_{\infty }} =\emptyset $. Therefore $(\d W\times
\bar V)\cap \Gamma_{f^{n_k}}=\emptyset $ for $k$ big enough.
This means that $\Phi_w\supset V$ as before and, moreover, $\Phi_w$ contains
the component of $C$ passing through $x$ minus the rest of $C$.

\smallskip To finish the proof let us distinguish two cases.

\smallskip\noindent{\slsf Case 1.} {\it Every component of $C$
contains such a point.} In this case our sequence $\{ f^{n_k}\} $
strongly converges on $X\setminus \{ \text{ finite set }\} $.
Furthermore, $\vol (\Gamma_{f^{n_k}})$ are uniformly bounded. Since
they can't be less than $(\deg f)^{n_k}\cdot \vol (X)$, we see that $f$ has
degree one, and $f_{\infty }$ is a degenerate map to $C$ because
$\Gamma_{f_{\infty}}$ has zero volume in this case. Moreover
$f_{\infty}$ cannot be holomorphic near $p$, otherwise $f^{n_k}$
would converge strongly in a neighborhood of $p$. $C$ in this case
should consist only from one component as a meromorphic image of an
irreducible variety.

\smallskip\noindent{\slsf Case 2.} {There exists a component $C_1$
of $C$ such that for all points $x\in C_1\setminus\cup_{i\not=
1}C_i$ $\dim_{(p,x)}\Gamma_{f_{\infty }}\cap (X\times \{x\}) >0$.
Then $f_{\infty }$ is a degenerate mapping of $X$ onto this $C_1$ and
therefore again $C_1$ is a single component of $C$. Indeed, 
any other component $C_2$ of $C$ should contain a point $x$ as above,
because the image of $f_{\infty}$ should be irreducible. I.e., in both
cases $C$ consists from one rational curve only.

\smallskip\qed

\smallskip The following simple example shows that the situation described in part
(b) of this Corollary can really happen. Let $X=\pp^2$ and $f:[z_0:z_1:z_2]\to [z_0:
2z_1:2z_2]$. Then for this $f$ we have the phenomena described above with
$p=[1:0:0]$ and  $C=\{ z_0=0\} $.

\newprg[FATOU.exmp]{Example}

Let us give one more example relevant to the Fatou sets.

\begin{exmp} \rm
\label{ex-deg-2}
Consider the following rational self-map of $\pp^2$:
\begin{equation}
\eqqno(map-f)
f : [z_0:z_1:z_2] \to [z_0^2z_1:z_1^{3}:z_0^2z_2].
\end{equation}
By induction one easily checks that
\begin{equation}
\eqqno(map-f-k)
f^k : [z_0:z_1:z_2] \to [z_0^{2^k}z_1^{2^k-1}:z_1^{2^{k+1}-1}:z_0^{2^{k+1}-2}z_2].
\end{equation}

\begin{figure}[h]
\centering
\includegraphics[width=2.5in]{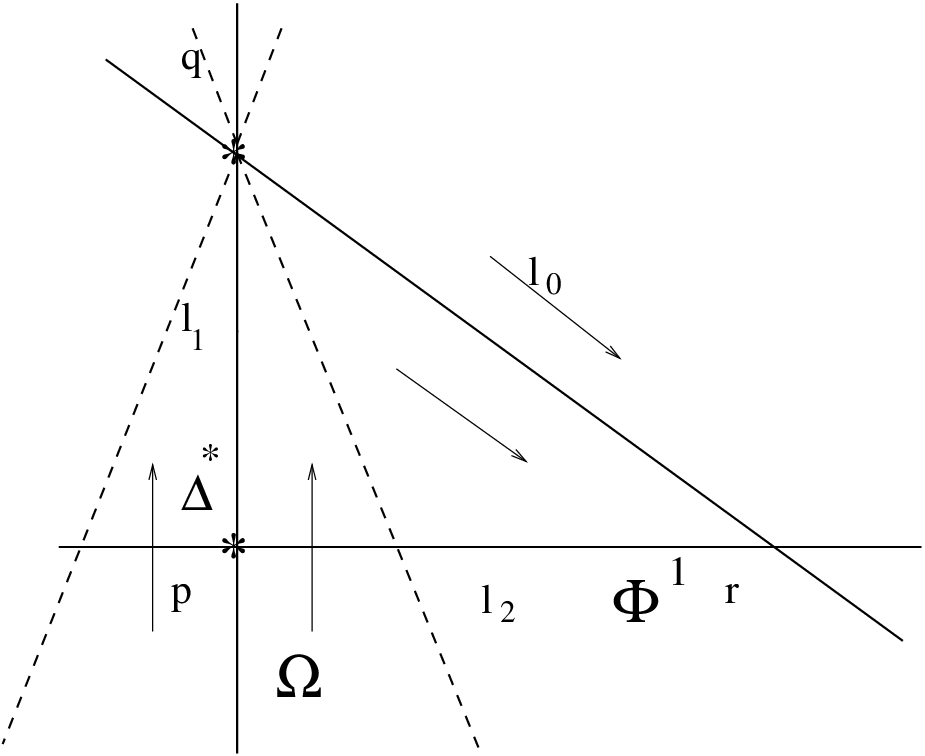}
\caption{Mapping $f$ contracts the line $l_1\deff \{z_1=0\}$ to the first of its points
of indeterminacy $q= [0:0:1]$, line at infinity $l_0\deff \{z_0=0\}$ to the regular point
$r=[0:1:0]$ and do not contracts anything to its second point of indeterminacy $p=[1:0:0]$.
Levi flat cone = Julia set for $f$ is marked by two punctured lines.} \label{pic1}
\end{figure}
Cover $\pp^2$ by three standard affine charts $U_i=\{z_i\not=0\}$ with coordinates
$u_1=\frac{z_1}{z_0}, u_2 = \frac{z_2}{z_0}$, $v_1 = z_0/z_1, v_2 = z_2/z_1$ and
$w_1=z_0/z_2, w_2 = z_1/z_2$ respectively. Mapping $f:U_0\to U_0$ writes as
\begin{equation}
\eqqno(map-f-u1)
f : (u_1,u_2) \to (u_1^2, u_2/u_1).
\end{equation}
We see from here that $f$ has degree $2$. Furthermore $f^k$ writes as
\begin{equation}
\eqqno(map-f-uk)
f^k : (u_1,u_2) \to (u_1^{2^k}, u_2/u_1^{2^k-1}).
\end{equation}
In the charts $f:U_1\to U_1$  our iterate writes as
\begin{equation}
\eqqno(map-f-vk)
f^k : (v_1,v_2) \to (v_1^{2^k}, v_1^{2^{k+1}-2}v_2) \to (0,0) = r \quad\text{on}\quad \{|v_1|<1\}.
\end{equation}
I.e., we see that $\Phi^1 = \{|v_1|<1\} = \{|u_1|>1\}$ is a component of the Fatou set of $f$, in all senses,
because all $f^k$ are holomorphic there and converge uniformly on compacts to a constant map to $r=
[0:1:0]$.

\smallskip Levi flat cone $L\deff \{|z_0|=|z_1|\}$ is a Julia set of $f$. It contains one of two indeterminacy
points of $f$, namely $q=[0:0:1]$. A connected component $\Omega$ of $\pp^2\setminus L$ different from $\Phi$
carries a more interesting information about $f^k$. First of all remark that as a mapping from $U_2$ to $U_2$
our iterate writes as
\begin{equation}
\eqqno(map-f-wk)
f^k : (w_1,w_2) \to \left(\frac{w_2^{2^k-1}}{w_1^{2^k-2}},\frac{w_2^{2^{k+1}-1}}{w_1^{2^{k+1}-2}}\right)
\to (0,0) = q \quad\text{on}\quad \{|w_2|<|w_1|\} = \{|u_1|<1, u_2\not=0\}.
\end{equation}
Therefore the second component $\Phi_s^2$ of the strong Fatou set contains the domain $\Omega\setminus
\{u_2=0\}$. Since it is easy from \eqqref(map-f-uk) to see that $f^k$ on compacts in $\Omega\setminus
\{u_2=0\}$ converge to $q$, and on the puncture disk $\Delta^*\deff \{u_2=0, 0<|u_1|<1\}$ to $p$, we
conclude that $\Phi_s^2 = \Phi_w^2=\Omega\setminus \{u_2=0\}$. Remark that the second component
$\Phi^2$ of $f$ in the sense of \cite{FS} is smaller, namely it is equal to $\Omega\setminus
\left(\{u_2=0\}\cup\{u_1=0\}\right)$, because the projective line $l_1\deff \{z_1=0\}$ is
the preimage of $I_f$ (and of all $I_{f^k}$) under $f$. Now let us turn to the second component
$\Phi_{\Gamma}^2$ of the $\Gamma$-Fatou set of $f$.

\begin{lem}
\label{gam-exmp}
For a fixed $0<\eps <1$ the volumes of graphs of $f_k$ over the bidisk $\Delta^2_{\eps}\subset U_0$
centered at $p$ are uniformly bounded. In particular $\Phi_{\Gamma}^2= \Omega$.
\end{lem}
\proof To estimate the volume of $\Gamma_{f_k}$ over a neighborhood
of $p$ we use coordinates $u_1,u_2$ and representation
\eqqref(map-f-uk). In these coordinates $\Delta^2_{\eps} =
\{u:\norm{u} <\eps\}$. Since $f^k$ preserves the vertical lines
$\{u_1=\const\}$ we can simplify our computations assuming that $f$
takes values in $\Delta\times \pp^1$, the last being equipped with
the Hermitian metric form $\omega = \omega_1 + \omega_2 =
\frac{i}{2}dz_1\wedge d\bar z_1 + \frac{i}{2}\frac{dz_2\wedge d\bar
z_2}{(1+|z_2|^2)^2}$. Now we get
\[
(f^k)^*\omega = \frac{i}{2}\left[2^{2k}|u_1|^{2^{k+1}-2} +
\frac{(1-2^k)^2|u_1|^{2^{k+1}-4}}{(|u_1|^{2^{k+1}-2} + |u_2|^2)^2}\right]du_1\wedge d\bar u_1
+ \frac{i}{2}|u_1|^{2^{k+1}-2} \frac{du_2\wedge d\bar u_2}{(|u_1|^{2^{k+1}-2} + |u_2|^2)^2}
\]
\[
+ \frac{i}{2}(1-2^k)\bar u_1^{-2^k}|u_1|^{2^{k+2}-4}\frac{du_2\wedge d\bar u_1}{(|u_1|^{2^{k+1}-2} +
|u_2|^2)^2} + \frac{i}{2}(1-2^k)u_1^{-2^k}|u_1|^{2^{k+2}-4}\frac{du_1\wedge d\bar u_2}{(|u_1|^{2^{k+1}-2} +
|u_2|^2)^2}.
\]
Therefore
\[
\int\limits_{\Delta^2_{\eps}}(f^k)^*\omega \wedge dd^c\norm{u}^2 =
\]
\[
=  4\pi^2\int\limits_0^{\eps}
\int\limits_0^{\eps}\left[2^{2k}r_1^{2^{k+1}-1} + (1-2^k)^2\frac{r_1^{2^{k+1}-3}}{(r_1^{2^{k+1}-2} + r_2^2)^2}
+ \frac{r_1^{2^{k+1}-1}}{(r_1^{2^{k+1}-2} + r_2^2)^2}\right]dr_1r_2dr_2 \le
\]
\[
\le \frac{4\pi^2\eps^{2^{k+1}}}{2^{2^{k+1}-2k}} + 2\pi^22^{2k}\int\limits_0^{\eps}\int\limits_0^{\eps^2}
\frac{r_1^{2^{k+1}-3}t}{(r_1^{2^{k+1}-2} + t)^2}dtdr_1 + 2\pi^22^{2k}\int\limits_0^{\eps}\int\limits_0^{\eps^2}
\frac{r_1^{2^{k+1}-1}}{(r_1^{2^{k+1}-2} + t)^2}dtdr_1 \le
\]
\[
\le 2\pi^22^{2k}\int\limits_0^{\eps}\int\limits_0^{\eps^2}\frac{r_1^{2^{k+1}-3}}{r_1^{2^{k+1}-2} + t}dtdr_1
+ 2\pi^2\int\limits_0^{\eps} \left(\frac{-1}{r_1^{2^{k+1}-2} + t}\Big|_0^{\eps^2}r_1^{2^{k+1}-1}\right)dr_1
\le
\]
\[
\le \pi^22^{3k+1}\int\limits_0^{\eps}r_1^{2^{k+1}-3}\ln\frac{1}{r_1}dr_1 + 2\pi^2\int\limits_0^{\eps}r_1dr_1.
\]
The second term is bounded and doesn't tend to zero as $k\to +\infty$. The first for $0<\eps <1$
can be obviously bounded by
\[
\pi^22^{3k+1}\int\limits_0^{\eps}r^{2^{k+1}-4}dr = \frac{\pi^22^{3k+1}}{2^{k+1}-3}\eps^{2^{k+1}-3}\to 0
\]
as $k\to +\infty$. And finally
\[
\int\limits_{\Delta^2_{\eps}}(f^k)^*\omega = \int\limits_{\Delta^2_{\eps}}\frac{2^{2k}|u_1|^{2^{k+2}-4}}
{(|u_1|^{2^{k+1}-2} + |u_2|^2)^2} + \frac{(2^{k}-1)^2|u_1|^{2^{k+2}-6}|u_2|^2}
{(|u_1|^{2^{k+1}-2} + |u_2|^2)^4} +
\]
\[
+ \frac{2(2^{k}-1)^2|u_1|^{2^{k+3}-8}|u_2|^2} {|u_1|^{2^{k+1}}(|u_1|^{2^{k+1}-2} + |u_2|^2)^4}d^4m(u) \le
\int\limits_{\Delta^2_{\eps}}\frac{2^{2k+2}|u_1|^{2^{k+2}-6}}{(|u_1|^{2^{k+1}-2} + |u_2|^2)^2}
d^4m(u)\approx
\]
\[
\approx 2^{2k+2} \int\limits_0^{\eps}\int\limits_0^{\eps^2}\frac{r_1^{2^{k+2}-5}}{(r_1^{2^{k+1}-2}+t)^2}dtdr_1
\le 2^{2k+2} \int\limits_0^{\eps}r^{2^{k+1}-3}dr = \frac{2^{2k+2}}{2^{k+1}-2}\eps^{2^{k+1}-2}\to 0
\]
as $k\to +\infty$. Therefore the lemma is proved.

\smallskip\qed

\smallskip Remark that the integral of $(f^k)^*\omega^2$ degenerates, as it should be, because $f^k$
$\Gamma$-converge on
$\Omega$ to a constant map. And to the contrary the integral of $(f^k)^*\omega$ doesn't
degenerate, moreover has order $\eps^2$. That means that bubbling takes place over all points of the
disk $\Delta^*=\{u_2=0, 0<|u_1|<1\}$. We see that for our map one has
\[
\Phi \subset \Phi_w=\Phi_s\subset \Phi_{\Gamma},
\]
and inclusions are strict.

\smallskip Finally let us see what is going on over the indeterminacy point $p=(0,0)$ in coordinates
$(u_1,u_2)$. Blowing up $2^{k}-1$ times at zero we see that $f^k$ stays degenerate on all exceptional
curves except the last one, which it send onto $l_1=\{u_1=0\}$. In appropriate coordinates $v_1,v_2$ on the
last blow up $f^k$ writes as
\[
\begin{cases}
u_1 = v_1^{2^k} \cr
u_2 = v_2.
\end{cases}
\]
Therefore the dynamical picture over $p$ can be described as follows. Let $\hat \Omega$ be an infinite blow up
of $\Omega$ over $p$ and let $C=\bigcup_{i=1}^{\infty}C_i$ be the Nori string of rational curves on $\hat\Omega$
over $p$. Then every $f^k$ lifts to a holomorphic map $\hat f^k:\hat \Omega \to \pp^2$ which is constantly equal
to $q$ on every $C_i$ except $C_{2^k-1}$. The last curve it sends bijectively onto the line $l_1$. At the same
time in the sense of divisors (currents) $f^k(C_{2^k-1}) = 2^kl_1$.

\end{exmp}

\begin{rema} \rm
\label{ex-deg-d}
{\bf a)} One can in the same way produce mappings of any given degree with the same properties
as in Example \ref{ex-deg-2}. It is sufficient to take
\begin{equation}
\eqqno(map-f-d)
f : [z_0:z_1:z_2] \to [z_0^dz_1:z_1^{d+1}:z_0^dz_2].
\end{equation}

\smallskip\noindent{\bf b)} Let us quote the
result of Maegava, see \cite{M}, which  shows that under an additional
assumption of "algebraic stability" of the dominant rational
self-map $f$ the Fatou set of \cite{FS} coincides with $\Phi_s$ and
$\Phi_w$.
\end{rema}

\def\entry#1#2#3#4\par{\bibitem[#1]{#1}
{\textsc{#2 }}{\sl{#3} }#4\par\vskip2pt}

\end{document}